\documentclass[a4paper, 
10pt
]{article}

\usepackage[a4paper,left=3cm,right=3cm,top=2.5cm,bottom=2.5cm]{geometry}
\usepackage{amsmath}
\usepackage{multirow}
\usepackage[table]{xcolor}
\usepackage{amsopn}
\DeclareMathOperator{\diag}{diag}
\DeclareMathOperator{\sech}{sech}

\usepackage{subcaption}
\usepackage{pgfplots}
\usepackage{pgfplotstable}
\usepackage{mathtools}
\usepackage{multicol}
\usepackage{comment}
\usepackage{booktabs}
\usepackage{amssymb}
\usepackage{url}
\usepackage{bm}
\usepackage[a4paper,left=3cm,right=3cm,top=2.5cm,bottom=2.5cm]{geometry}
\usepackage{amsmath}
\usepackage{multirow}
\usepackage{float}
\usepackage{tabularx}
\usepackage{enumerate}
\usepackage{array}
\usepackage{xspace}
\usepackage{tikz}
\usepackage[export]{adjustbox}
\usepackage{tikzsymbols}
\usepackage{authblk}

\numberwithin{equation}{section}
\usepgfplotslibrary{groupplots}
\usepgfplotslibrary{fillbetween}
\pgfplotsset{compat=1.5}
\usetikzlibrary{calc,trees,positioning,arrows,chains,shapes.geometric,%
    decorations.pathreplacing,decorations.pathmorphing,shapes,%
    matrix,shapes.symbols, decorations.markings, patterns,fit,quotes,%
    arrows.meta, spy}

\usepackage{xcolor}
\definecolor{Gray}{gray}{0.9}

\usepackage{hyperref}
\usepackage{siunitx}
\usepackage{cleveref}
\usepackage{algorithmic}
\usepackage{algorithm}
\usepackage[bottom]{footmisc}

\graphicspath{{figures/}}

\newcommand{\mupar}{\ensuremath{\boldsymbol{\mu}}}

\newcommand{\ntrain}{p}
\newcommand{\timedim}{N}
\newcommand{\spacedim}{m}
\newcommand{\npod}{n}

\newcommand{\RA}[1]{{\color{black}#1}}
\newcommand{\RB}[1]{{\color{black}#1}}

\newcommand{\matr}[1]{\mathbf{#1}}

\begin{document}

\title{A dynamic mode decomposition extension for the forecasting of parametric dynamical systems} 

\author[]{Francesco Andreuzzi\footnote{francesco.andreuzzi@sissa.it}}
\author[]{Nicola Demo\footnote{nicola.demo@sissa.it}}
\author[]{Gianluigi Rozza\footnote{gianluigi.rozza@sissa.it}}

\affil{Mathematics Area, mathLab, SISSA, via Bonomea 265, I-34136
  Trieste, Italy}

\maketitle

\begin{abstract}
  Dynamic mode decomposition (DMD) has recently become a popular tool for the non-intrusive analysis of dynamical systems. Exploiting the proper orthogonal decomposition as dimensionality reduction technique, DMD is able to approximate a dynamical system as a sum of (spatial) basis evolving linearly in time, allowing for a better understanding of the physical phenomena or for a future forecasting.
  We propose in this contribution an extension of the DMD to parametrized dynamical systems, focusing on the future forecasting of the output of interest in a parametric context. Initially, all the snapshots --- for different parameters and different time instants --- are projected to the reduced space, employing the DMD (or one of its variants) to approximate the reduced snapshots for a future instants. Still exploiting the low dimension of the reduced space, the predicted reduced snapshots are then combined using a regression technique, enabling the possibility to approximate any untested parametric configuration in any future instant.
  We are going to present here the algorithmic core of the aforementioned method, presenting at the end three different test cases with incremental complexity: a simple dynamical system with a linear parameter dependency, a heat problem with nonlinear parameter dependency and a fluid dynamics problem with nonlinear parameter dependency.
\end{abstract}


\section{Introduction}
\label{sec:intro}
Several applications in the \RB{field of} computational sciences require an iterative or real--time evaluation of the mathematical model in hand. In these contexts, Reduced Order Modeling (ROM) \RB{constitutes a popular tool for boosting the} performance of the methods by reducing the dimensionality of the original model. Among all ROM methods, Dynamic Mode Decomposition (DMD) \cite{dmd} has lately gained popularity for \RB{what concerns the simulation} of dynamical systems, \RB{thanks to the wide applicability of the method in multiple contexts, e.g.} fluid dynamics problems \cite{DemoOrtaliGustinRozzaLavini2020BUMI, TezzeleDemoMolaRozza2018, RozzaMalikDemoTezzeleGirfoglioStabileMola2018ECCOMAS, gadalla2021les}. \RB{DMD provides dimensionality reduction} in a data-driven fashion, which then enables locating dominant structures in the original problem and finally predicting the system state in future \RB{time instants}.

However, in many applications, it is fundamental to study not only the dynamics of the system but also its parametric dependency. For practical examples, the reader can think of the position of the source in an unsteady heat problem, or the fluid velocity in a computational fluid dynamics one. \RB{The task, in this case, is usually to predict the evolution of the system for both unseen parameters and future time instants}.

There have been several attempts to achieve this goal, which we briefly summarize in this paragraph. Parametric time-dependent problems can be treated in a projection framework exploiting Proper Orthogonal Decomposition (POD)~\cite{georgaka2018parametric,hijazi2020data} for dimensionality reduction. Such an approach requires --- due to the projection step --- the knowledge of the original model, as well as access to the discrete operators characterizing the problem. \RB{Earlier attempts were made with non-intrusive methodologies as well. Extracting the main spatio-temporal structure of the problem to isolate the parameter dependency was proposed in~\cite{hoang2021projection}. Even though this technique is not able to predict the output of the system for future time steps, it enables the reconstruction of the parametric manifold. A different direction was probed in \cite{brunton2014compressive}, where authors employ POD to classify the dynamic region of the current sample and provide a low-rank approximation by means of Galerkin projection. A second attempt in this direction is presented in \cite{kramer2017sparse}: the authors employed DMD to analyze and classify complex flow in different regimes, thus providing a coarse reconstruction of the system's state by means of a low-order flow approximation. A pipeline similar to ours is depicted in~\cite{amsallem2011online}, where authors employ projection-based ROM and interpolation on a reduced space. However, they do not take into account DMD or any other time forecasting techniques. Finally, a parametric version of DMD has been already studied in~\cite{pdmd}, however, the aim of that work is the extraction of coherent structures shared between the different configurations of the system, rather than obtaining predictions on unseen parameters in future time instants}.

In this work, we propose an extension of a \RB{previous work presented in~\cite{TezzeleDemoStabileMolaRozza2020MOR}, which consists in the application of DMD to parameter-dependent scalar output in order to predict such output for untested parameters and time steps. Our aim is to predict vectorial fields as well, thus achieving the general goal task outlined above}. To the best of the authors' knowledge, an approach for predicting the fields of interest for unseen time instants and parameters is not explored in the literature yet.

Following the standard DMD method, our extension takes in input an initial dataset of high--fidelity snapshots, equispaced in time and for different parametric configurations. Such snapshots are then reduced in dimensionality by means of POD, and then re-arranged in order to feed the DMD algorithm. DMD is applied in order to forecast the reduced snapshots in the future. This step output reduced snapshots for all the parametric samples which compose the initial dataset. Then, the parametric manifold can be reconstructed using an approach similar to POD with interpolation \cite{farhat, DemoOrtaliGustinRozzaLavini2020BUMI} by means of a regression technique whose aim is approximating the mapping between parameters and reduced snapshots. Since this technique works on a reduced space, the computational cost for any new parametric evaluation \RB{(a popular use case in the field)} is small with respect to traditional discretization methods like Finite Elements and Finite Volumes. Moreover, the method is equation-agnostic, since it relies only on the high--fidelity snapshots, thus not requiring any kind of information about the original model \RB{(except for the parameters that generated the observed samples)}.
\RA{The present work will focus also on a novel approach that exploits the knowledge of the problem at hand in order to stabilize the DMD operator, such that forecasting in the future --- the goal of this work --- exhibits a greater accuracy with respect to standard operator.}

The present contribution is organized as follows: in \cref{sec:methods} \RB{we briefly introduce DMD and a few of the variants taken into account in our experiments}, and in \cref{meth:paramdmd} \RB{we discuss in detail our novel pipeline for extending DMD to parametric time-dependent problems}. \RA{\Cref{meth:stab} is devoted to the discussion of the stabilization approach applied to the DMD operator, while} in \cref{sec:results} we present numerical results on three problems with heterogeneous parametric dependency. Finally, in \cref{sec:conclusions} we discuss several possible future developments \RB{for subsequent works}.

\section{Methods}
\label{sec:methods}
\label{sec:dmd}
We dedicate this section to a brief introduction to Dynamic Mode Decomposition \RA{and some of its variants already introduced in the literature, in order of providing an algorithmic overview to the reader. The original algorithm, even if slightly modified, is indeed the backbone of the parametric extension here presented. We specify moreover the capability to exploit the introduced approach also with DMD variants, extending the possibility of applications.}\\

DMD is a data-driven method for the analysis of dynamic systems. We provide a brief introduction to the general technique in order to make the extension described in \cref{sec:methods} more clear to the reader. We consider a time-dependent vector function $\RA{x(t)} \in \mathbb{R}^\spacedim$ and we call $\matr{x}_{\RA{t^k}}$ the \emph{snapshot} of the vector quantity of interest collected ad time instant $t^k$.
We seek an operator $\matr{A} \in \mathbb{R}^{\spacedim \times \spacedim}$ such that:
\begin{gather}\label{intro:dmd_def}
    \matr{x}_{t^{k+1}} = \matr{A} \matr{x}_{t^k},
\end{gather}
\RA{where $\matr{x}_{t^k}$, $\matr{x}_{t^{k+1}}$ are the two snapshots collected in equential time instants. 
}
The reader could think of the operator $\matr{A}$ as a finite approximation of the infinite-dimensional \emph{Koopman operator}~\cite{koopman}. In order to build such an operator satisfying \cref{intro:dmd_def} up to an acceptable degree of approximation, we consider a set $\mathcal{T}$ of discrete time instants equispaced in time, which we are going to refer to as \emph{tested} -- or \emph{known} -- time instants. These equispaced time instants are labeled with cardinal numbers $\{1, \dots, \timedim\}$ (therefore we assume that the cardinality of $\mathcal{T}$ is $\timedim$). Thus we take $\timedim$ snapshots $\{\matr{x}_{\RA{t^1}}, \dots, \matr{x}_{\RA{t^\timedim}}\}$ of the function of interest, to be used during the training phase of DMD. We arrange the snapshots in two matrices $\matr{X},\matr{Y} \in \mathbb{R}^{\spacedim \times \timedim - 1}$ defined as follows:
\begin{gather}
    \matr{X} = \begin{bmatrix}
        \matr{x}_{\RA{t^1}} & \matr{x}_{\RA{t^2}} & \dots & \matr{x}_{\RA{t^{\timedim -1}}}
    \end{bmatrix}, \qquad
    \matr{Y} = \begin{bmatrix}
        \matr{x}_{\RA{t^2}} & \matr{x}_{\RA{t^3}} & \dots & \matr{x}_{\RA{t^\timedim}}
    \end{bmatrix}.
\end{gather}
We now attempt to find the matrix $\matr{A}$ such that $\matr{Y} \approx \matr{A}\matr{X}$ (note that such an operator satisfies the condition \cref{intro:dmd_def} $\forall k \in \{1,\dots,\timedim-1\}$). Denoting the Moore-Penrose pseudo-inverse operator by $^\dagger$, we can find the best-fit operator $\matr{A} = \matr{Y}\matr{X}^{\dagger}$ \cite{penrose1956best}. However this computation is expensive and generally intractable, therefore we follow another path.

Consider the $r$-ranked singular value decomposition of $\matr{X} \approx \matr{U}_r \matr{\Sigma}_r \matr{V}_r^*$ \cite{quarteroni}. Since the matrix $\matr{U}_r$ is orthogonal, we take the low-rank projection $\matr{\widetilde{A}} \in \mathbb{R}^{r \times r}$ of the operator $\matr{A}$ onto the subspace spanned by the columns of $\matr{U}_r$:
\begin{gather}
    \matr{\widetilde{A}} = \matr{U}_r^* \matr{A} \matr{U}_r = \matr{U}_r^* \matr{Y} \matr{X}^{\dagger} \matr{U}_r
    = \matr{U}_r^* \matr{Y} \matr{V}_r \matr{\Sigma}_r^{-1} \matr{U}_r^* \matr{U}_r = \matr{U}_r^* \matr{Y} \matr{V}_r \matr{\Sigma}_r^{-1}.
    \label{eq:dmd_fit}
\end{gather}
Note that this computation does not involve the evaluation of the unprojected operator $\matr{A}$. It can be shown that the eigenvalues of the full-dimensional operator $\matr{A}$ are the same of the low-dimensional operator $\matr{\widetilde{A}}$, and that we can use the eigenvectors of the former to compute the eigenvectors of the latter \cite{dmd}. Since the eigendecomposition of the low-rank operator involves a tractable computation due to the low dimensionality of the matrix $\matr{\widetilde{A}}$, the eigendecomposition of the full-dimensional operator $\matr{A}$ is easily obtainable through the manipulation of low-rank matrices.

Assuming that the eigendecomposition of the full-dimensional operator $\matr{A}$ yields the matrix $\matr{\Phi} \in \mathbb{R}^{\spacedim \times r}$ of eigenvectors of $\matr{A}$, and the diagonal matrix $\matr{\Lambda} \in \mathbb{R}^{r \times r}$ of its eigenvalues, it can be shown that:
\begin{gather} \label{eq:dmd_future}
    \matr{x}_{\RA{t^k}} = \matr{\Phi} \matr{\Lambda}^k \matr{\Phi}^\dagger \matr{x}_{\RA{t^1}}.
\end{gather}
There are virtually no constraints on $k \in \mathbb{N}$, therefore it is possible to use \cref{eq:dmd_future} to extrapolate the value of $\matr{x}_{\RA{t^k}}$ for $k > \timedim$, namely we can predict the behavior of the function of interest outside the tested time window, clearly, with some uncertainty~\cite{lu2020prediction}.

\subsection{Higher-Order DMD variant}
\label{sec:hodmd}
The quality of results obtained using DMD decays if $\spacedim \approx \timedim$ \cite{hodmd}, where $\spacedim$ is the space dimension (the number of components of the vector function of interest) and $\timedim$ is the number of time instants considered during the training phase. It has been observed that the situation becomes even worse when $\spacedim < \timedim$. A possible cure for this problem is provided by Higher Order variants of DMD \cite{hodmd, vasconcelos2019dynamic}.

The key part of those variants is the replacement of the usual expression $x_{\RA{t^{k+1}}} \approx \mathbf A x_{\RA{t^{k}}}$ with the following expression, which uses $d$ time lagged snapshots:
\begin{gather}\label{eq:hodmd}
\matr{x}_{\RA{t^{k+d}}} = \matr{A}_1 \matr{x}_{\RA{t^{k}}} + \matr{A}_2 \matr{x}_{\RA{t^{k+1}}} + \dots + \matr{A}_d \matr{x}_{\RA{t^{k+d-1}}}.
\end{gather}
This is usually called \emph{higher order Koopman assumption}. Then we set:
\begin{gather}
    \bar{\matr{x}}_{\RA{t^k}} \coloneqq \begin{bmatrix} \matr{x}_{\RA{t^k}} & \matr{x}_{\RA{t^{k+1}}} & \dots & \matr{x}_{\RA{t^{k+d-1}}}\end{bmatrix}^\intercal
\end{gather}
This enables rewriting \cref{eq:hodmd} in a compact form which resembles \cref{intro:dmd_def}: indeed it is clear that $\bar{\matr{x}}_{\RA{t^{k+1}}} = \matr{A} \bar{\matr{x}}_{\RA{t^k}}$. In fact this form enables finding an approximation of the operator $\matr{A}$ like we did in \cref{sec:dmd}. Experimental results show that this approach produces better results in several settings \cite{hodmdnoisy,hodmdflow}.

Higher Order DMD variants are useful additions to our toolchain since we usually consider a significant dimensionality reduction before the application of DMD. This is a compulsory step because otherwise, the interpolation would be unfeasible from a computational perspective otherwise.
\RB{
\subsection{Noise-tolerant DMD variants}
\label{meth:noise}
DMD results could be damaged by the noisy data~\cite{bagheri2013effects}, even if the reduced space is used for projecting the DMD operator. Moreover, in our pipeline, we consider multiple steps which could propagate and in some cases increase the effect of noise. Along with the simple methodology depicted in \cref{meth:stab}, a feasible way to improve the quality of results in presence of sample noise is using noise-tolerant DMD variants. The implementation style of our parametric DMD variant in \cite{pydmd} enables easy selection of several DMD variants, thus easing fine-tuning the DMD variant choice and the corresponding hyperparameters.

In this brief paragraph, we briefly mention a few works which could improve the results of our algorithm. In \cite{dawson2016characterizing} the authors propose a few techniques to improve DMD results in presence of noise, by means of well-known methods in sound processing, combining forward and backward in time DMD, and via a novel total-least-squares DMD algorithm. Total DMD~\cite{hemati2017biasing} claims to be able to incorporate noise information to \emph{de-bias} the usual subspace projection step in DMD by means of an augmented snapshot matrix. Finally, statistical \emph{bagging} is proposed as an alternative solution for noise. This way of dealing with the problem is explored in~\cite{sashidhar2022bagging}, whose authors propose the Bagging OPtimized DMD (BOP-DMD). It's important to spend a few words also on time measurement noise since problematic datasets are not always caused by bad sensors or uncertain measurements. An additional source of problems is given by measurements performed in uneven time instants. This issue is addressed in~\cite{askham2018variable}, where the authors pair the DMD with the variable projection method for a nonlinear least squares problem. We decided not to include an in-depth analysis of noise-tolerant DMD variants to avoid overloading the current work. We reserve to pursue this direction further in future developments.
}
\section{Dynamic Mode Decomposition for parametric problems}
\label{meth:paramdmd}

In this section we propose a possible application of DMD to a parametrized dynamical system. The parameter belongs to a generic set $\mathbb{P}$, which we refer to as \emph{parameters set}. Considering a discrete space, we define the output of interest as the vector function $\mathbf{x}^\mu_t \in \mathbb{R}^\spacedim$, which depends on the parameter $\mu \in \mathbb{P}$ and on the time instant $t$.

As in the original formulation, we initially collect a set of $\timedim$ high-fidelity snapshots of the system with a fixed parameter at equispaced time instants. We recall that $\mathcal{T}$ is the set of time instants used during the training phase, we keep using the labeling of time instants introduced in \cref{sec:dmd}, namely $\mathcal{T} = \{1, \dots, \timedim\}$. We also consider $\ntrain$ parameters, which we collect into the \emph{training set} $\mathcal{S} = \{\mu_1, \dots, \mu_\ntrain\}$, and we replicate the previous step $\ntrain$ times, one for each parametric sample. We obtain then $\ntrain$ parametric configurations, each one containing $\timedim$ time instants, overall $\ntrain \times \timedim$ snapshots. The basic idea is then to apply a POD-based reduction to the snapshots in order to obtain their reduced representation, also referred to in the literature as modal coefficients or reduced snapshots \cite{rozza_pod}. Then we apply DMD in order to approximate reduced snapshots in future time instants \cite{dmd}: this allows us to explore the solution manifold using some regression techniques \cite{farhat} in future time instants, approximating the relation between parameters and reduced snapshots. Finally, the interpolated reduced snapshot is mapped back to the full-dimensional space.

\begin{figure}[h]
    \centering
    \includegraphics[width=\textwidth]{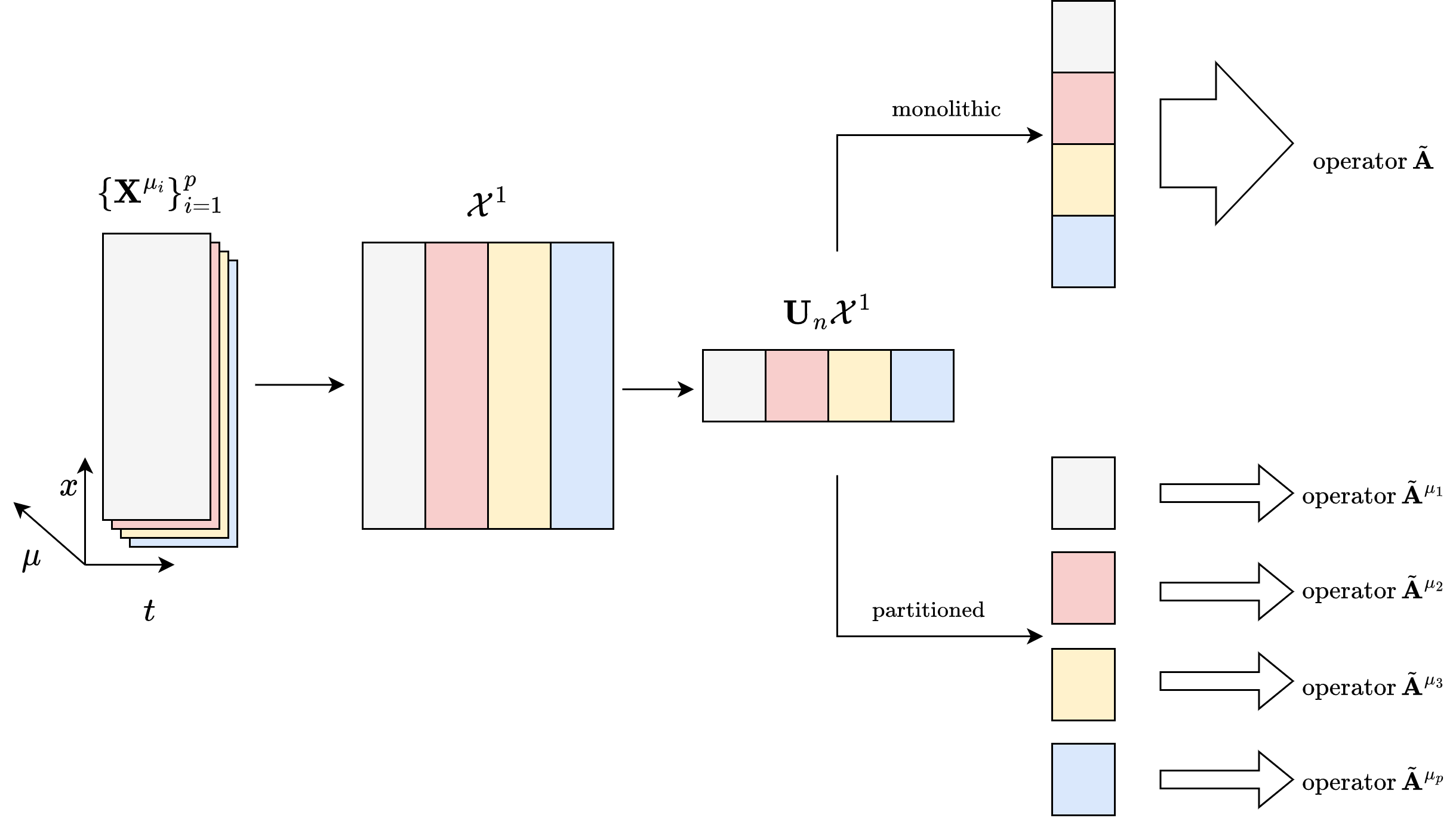}
    \caption{\RA{Scheme of the numerical pipeline to perform the offline stage, isolating the differences between \emph{monolithic} and \emph{partitioned}.}}
    \label{fig:offline}
\end{figure}

During the numerical experiments, we consider two different approaches for the application of DMD to the reduced snapshots. In the first case, we build a unique linear operator to fit the dynamics of the entire (parameterized) system. Alternatively, $\ntrain$ linear operators are constructed to approximate the dynamics of the $\ntrain$ parametric configurations. \RA{
Figure~\ref{fig:offline} visually schematized these two alternatives, highlighting similarities and differences between them.
}
In the next subsections, we formally define these two approaches, referred from now on as \emph{monolithic} and \emph{partitioned}, presenting also the regression technique and some numerical stabilization we impose on DMD operator for more robust forecasting.

To conclude this preliminary section, we highlight that both the versions we are going to explain have been implemented in the PyDMD package~\cite{pydmd}, an open-source Python package collecting several implementations of DMD and its variants.
\subsection{Monolithic approach}
\label{meth:full}

We consider $\ntrain$ matrices $\matr{X}^{\mu_i}$ (one for each parameter in the training set $\mathcal{S} = \{\mu_1, \dots, \mu_\ntrain\}$) and the matrix $\matr{\mathcal{X}}_1$ defined as follows:
\begin{gather}\label{eq:full_arrange_snapthots}
    \matr{X}^{\mu_i} \coloneqq \begin{bmatrix}
        \mid & \dots & \mid\\
        \matr{x}_{\RA{t^1}}^{\mu_i} & \dots & \matr{x}_{\RA{t^\timedim}}^{\mu_i}\\
        \mid & \dots & \mid\\
    \end{bmatrix} \in \mathbb{R}^{\spacedim \times \timedim}, \quad \matr{\mathcal{X}}_1 \coloneqq \begin{bmatrix}
        \matr{X}^{\mu_1} \dots \matr{X}^{\mu_\ntrain}
    \end{bmatrix} \in \mathbb{R}^{\spacedim \times \timedim \ntrain}.
\end{gather}
Note that the columns of each $\matr{X}^{\mu_i}$ are the (sampled) evolution over time of the vector function $\matr{x}_{\RA{t^k}}^{\mu_i} \in \mathbb{R}^{\spacedim}$, in the time instants belonging the set $\mathcal{T}$.
We apply POD to the matrix $\matr{\mathcal{X}}_1$ to obtain the optimal\footnote{in linear sense} space of dimension $\npod$ to represent all the original snapshots, with $ 0 < \npod \leq \min(\spacedim, \timedim \ntrain)$. POD method returns indeed the matrix $\matr{U}_\npod \in \mathbb{R}^{\spacedim \times \npod}$ which contains the reduced space basis, or POD modes~\cite{rozza_pod}. The number $\npod$ of retained POD modes is selected by the user and must be selected carefully since it affects the quality of the results provided by the algorithm. A small value may result in a low-quality reconstruction of the output function, since important spacial information may have been discarded; a high value may increase the computational cost of the online phase, as well as increase the error introduced by regressors due to the higher dimension.

We can exploit the matrix $\matr{U}_\npod$ to compute the modal coefficients corresponding to the snapshots stored in the columns of the matrix $\matr{\mathcal{X}}_1$: in matrix form, we have $\widetilde{\matr{X}} = \matr{U}_\npod^* \matr{\mathcal{X}}_1 \in \mathbb{R}^{\npod \times \timedim \ntrain}$. The columns of $\widetilde{\matr{X}}$ are the reduced snapshots that we mentioned earlier. We refer by $\widetilde{\matr{x}}_{\RA{t^k}}^{\mu_i} \coloneqq \mathbf{U}_\npod^*\matr{x}_{\RA{t^k}}^{\mu_i} \in \mathbb{R}^{\npod}$ to the vector of POD coefficients corresponding to the $k$-th snapshot of the system with parameter $\mu_i$. Note that POD coefficients evolve in time, therefore we would like to predict this evolution in future time instants in order to obtain the full-dimensional counterpart of the snapshot, which we could later use to explore the solution manifold.

In order to employ the DMD method introduced in \cref{sec:dmd} we arrange POD coefficients in the following way:
\begin{gather}\label{eq:pod_coefficients_matrix}
    \widetilde{\matr{X}}^{\mu_i} \coloneqq \mathbf{U}_\npod^\intercal \matr{X}^{\mu_i} = \begin{bmatrix}
        \mid & \dots & \mid\\
        \widetilde{\matr{x}}_{\RA{t^1}}^{\mu_i} & \dots & \widetilde{\matr{x}}_{\RA{t^\timedim}}^{\mu_i}\\
        \mid & \dots & \mid\\
    \end{bmatrix} \in \mathbb{R}^{\npod \times \timedim}, \qquad \mathcal{X}_2 \coloneqq \begin{bmatrix}
        \widetilde{\matr{X}}^{\mu_1}\\
        \widetilde{\matr{X}}^{\mu_2}\\
        \vdots\\
        \widetilde{\matr{X}}^{\mu_\ntrain}
    \end{bmatrix} \in \mathbb{R}^{\npod \ntrain \times \timedim}.
\end{gather}
Each matrix $\widetilde{\matr{X}}^{\mu_i}$ is the reduced representation (in the POD subspace) of the corresponding matrix $\matr{X}^{\mu_i}$ and is obtained with a single matrix multiplication with POD modes. At this point, we apply DMD to the matrix $\matr{\mathcal{X}}_2$, whose columns are temporal snapshots representing different realizations of the parametric system. Arranging these columns in two matrices, $\matr{\mathcal{X}}_2^1 \in \mathbb{R}^{\npod \ntrain \times (\timedim-1)}$ and $\matr{\mathcal{X}}_2^2 \in \mathbb{R}^{\npod \ntrain \times (\timedim-1)}$, we recall that the DMD operator is the matrix $\widetilde{\matr{A}}$ such that $\matr{\mathcal{X}}_2^2 = \widetilde{\mathbf{A}}\matr{\mathcal{X}}_2^1$ as we stated in \cref{sec:dmd}. We note that, since in our method the dimension reduction is employed at the beginning, the typical compression step performed in the standard DMD formulation (eq.~\eqref{eq:dmd_fit}) might not be needed, depending by the number of parametric samples. In this case, the DMD operator can be found using just the least-square procedure.

As we mentioned the DMD operator $\widetilde{\matr{A}}$ can be used to advance a vector of reduced snapshots by one time instant:
\begin{gather*}
    \begin{bmatrix}
        \widetilde{\matr{x}}_{\RA{t^{k+1}}}^{\mu_1}\\
        \vdots\\
        \widetilde{\matr{x}}_{\RA{t^{k+1}}}^{\mu_\ntrain}
    \end{bmatrix}
    =
    \widetilde{\matr{A}} \begin{bmatrix}
        \widetilde{\matr{x}}_{\RA{t^k}}^{\mu_1}\\
        \vdots\\
        \widetilde{\matr{x}}_{\RA{t^k}}^{\mu_\ntrain}
    \end{bmatrix} \in \mathbb{R}^{\npod \ntrain}.
\end{gather*}
It might not be formally correct to keep using the notation $\widetilde{\matr{x}}_{\RA{t^k}}^{\mu_i}$ when $k > \timedim$, since in this case, the vector is an approximation of the reduced snapshot computed using DMD and \cref{eq:dmd_future}. However, we are going to use the same notation to avoid introducing another one for predicted reduced snapshots. We do this without ambiguity since in this work we are not going to consider reduced snapshots for $k > \timedim$.

We rewrite the last equation using the notation adopted, and we also replace the generic $k$ with the first time instant of interest for our discussion, namely the last one:
\begin{gather}\label{eq:predict_future_pod}
    \begin{bmatrix}
        \widetilde{\matr{x}}_{\RA{t^{\timedim+1}}}^{\mu_1}\\
        \vdots\\
        \widetilde{\matr{x}}_{\RA{t^{\timedim+1}}}^{\mu_\ntrain}
    \end{bmatrix}
    =
    \widetilde{\matr{A}} \begin{bmatrix}
        \widetilde{\matr{x}}_{\RA{t^\timedim}}^{\mu_1}\\
        \vdots\\
        \widetilde{\matr{x}}_{\RA{t^\timedim}}^{\mu_\ntrain}
    \end{bmatrix} \in \mathbb{R}^{\npod \ntrain}.
\end{gather}

\begin{algorithm}
    \caption{Dynamic Mode Decomposition for parametric problems | Monolithic }
    \hspace*{\algorithmicindent} \textbf{Input}: $\matr{X}^{\mu_1}, \dots, \matr{X}^{\mu_\ntrain} \in \mathbb{R}^{\spacedim \times \timedim}$\\
    \label{alg:pdmd_monolithic}
    \begin{algorithmic}
        \STATE{Construct the matrix $\matr{\mathcal{X}}_1 \in \mathbb{R}^{\spacedim \times \timedim \ntrain}$}
        \STATE{Compute the POD coefficients of the columns of $\matr{\mathcal{X}}_1$}
        \STATE{Store the coefficients into $\widetilde{\matr{X}}^{\mu_i}, i \in \{1, \dots, \ntrain\}$}
        \STATE{Construct the matrix $\matr{\mathcal{X}}_2 \in \mathbb{R}^{\npod \ntrain \times \timedim}$}
        \STATE {Compute the DMD operator $\widetilde{\matr{A}}$ by fitting $\matr{\mathcal{X}}_2$}
    \end{algorithmic}
\end{algorithm}
\noindent We remark that the monolithic approach relies on a unique DMD operator to express the dynamics of the parametric system. This implies that the operator is able to detect recurrent patterns in the dynamics of different parametric configurations, making its usage \RA{particularly} profitable when the different realizations share common behaviors.

\subsection{Partitioned approach}
\label{meth:partial}
We considered also another similar approach, inspired by the observation of patterns extracted with DMD as shown in \cref{meth:full} for particular dynamical systems. We are going to use the same notation, when possible, introduced in \cref{meth:full}.

Starting from the $\ntrain$ matrices $\matr{X}^{\mu_i}$ we construct $\matr{\mathcal{X}}_1$ as in \cref{eq:full_arrange_snapthots}. Then we apply POD to $\matr{\mathcal{X}}_1$ to obtain POD modes $\matr{U}_\npod \in \mathbb{R}^{\spacedim \times \npod}$, and we compute modal coefficients $\widetilde{\matr{X}}^{\mu_i} = \mathbf{U}_\npod^\intercal \matr{X}^{\mu_i}$ of the full order snapshots, for each parameter in the training set $\mathcal{S} = \{\mu_1, \dots, \mu_\ntrain\}$. Note that this procedure is the same which we presented in \cref{meth:full}. However, at this point we do not build the matrix $\matr{\mathcal{X}}_2$: instead we perform $\ntrain$ separate DMDs, one for each matrix $\widetilde{\matr{X}}^{\mu_i}, i \in \{1,\dots,\ntrain\}$. Namely, for each $\mu_i$ with $i \in \{1,\dots,\ntrain\}$ we consider the matrices $\widetilde{\matr{X}}^{\mu_i}_1, \widetilde{\matr{X}}^{\mu_i}_2 \in \mathbb{R}^{\npod \times (\timedim-1)}$ which results from taking $\timedim-1$ columns of $\widetilde{\matr{X}}^{\mu_i}$ --- as in the standard DMD method ---  and we approximate the DMD operator $\widetilde{\matr{A}}^{\mu_i} \in \mathbb{R}^{\npod \times \npod}$ such that $\widetilde{\matr{X}}^{\mu_i}_2 = \widetilde{\matr{A}}^{\mu_i} \widetilde{\matr{X}}^{\mu_i}_1$. After that the prediction of future POD coefficients is carried out independently for each tested parameter using \cref{eq:predict_future_pod}.

\begin{algorithm}
    \caption{Dynamic Mode Decomposition for parametric problems | Partitioned}
    \hspace*{\algorithmicindent} \textbf{Input}: $\matr{X}^{\mu_1}, \dots, \matr{X}^{\mu_\ntrain} \in \mathbb{R}^{\spacedim \times \timedim}$\\
    \label{alg:pdmd_partitioned}
    \begin{algorithmic}
        \STATE{Construct the matrix $\matr{\mathcal{X}}_1 \in \mathbb{R}^{\spacedim \times \timedim \ntrain}$}
        \STATE{Compute the POD coefficients of the columns of $\matr{\mathcal{X}}_1$}
        \STATE{Store the coefficients into $\widetilde{\matr{X}}^{\mu_i}, i \in \{1, \dots, \ntrain\}$}
        \FOR{$i \in \{1, \dots, p\}$ }
          \STATE {Compute the DMD operator $\widetilde{\matr{A}}^{\mu_i}$ by fitting $\widetilde{\matr{X}}^{\mu_i}$}
        \ENDFOR
    \end{algorithmic}
\end{algorithm}
\noindent According to our experiments this approach is more reliable than that presented in \cref{meth:full} when applied on systems for which the variation of the parameter leads to a considerable change in the behavior. Different instances of these systems (i.e. realizations with different parameters) are often weakly coupled, in the sense that it is usually hard to find recurrent patterns, if not possible at all. The application of the method presented in this section partially fixes the problem.

\subsection{Online phase}
Using the features presented in \cref{sec:dmd} and the methods proposed in \cref{meth:full} and \cref{meth:partial} we are able to exploit \cref{eq:dmd_future} to predict the value of POD coefficients in future time instants. To simplify the discussion we are going to assume that the approach is used to predict just one time instant, namely $t^*$. It will be clear that this discussion can be easily extended to a more general case.

After the \emph{offline phase} we can exploit the DMD operator (or the several DMD operators in case of the partitioned approach) to predict the reduced snapshots at time $t^*$ for the parameters in the training set $\mathcal{S} = \{\mu_1, \dots, \mu_\ntrain\}$, obtaining then the matrix $\matr{\mathcal{X}}_3$ which has the following form:
\begin{gather}
    \matr{\mathcal{X}}_3 = \begin{bmatrix}
        \mid & \dots & \mid\\
        \widetilde{\matr{x}}_{\RA{t^*}}^{\mu_1} & \dots & \widetilde{\matr{x}}_{\RA{t^*}}^{\mu_\ntrain}\\
        \mid & \dots & \mid\\
    \end{bmatrix} \in \mathbb{R}^{\npod \times \ntrain}
\end{gather}
The matrix $\matr{\mathcal{X}}_3$ has one row for each retained POD coefficient and one column for each realization of the parametric system used to train the algorithm.
\begin{figure}[h]
    \centering
    \includegraphics[width=\textwidth]{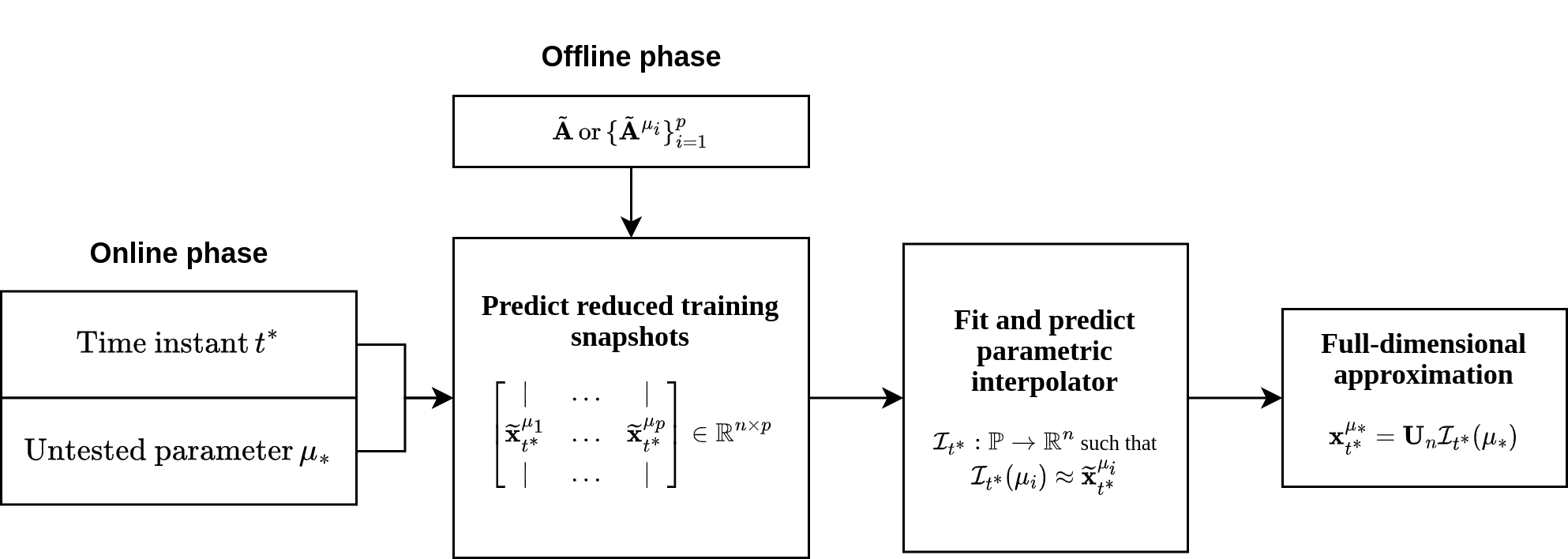}
    \caption{\RA{Scheme of the numerical pipeline to perform the online stage.}}
    \label{fig:online}
\end{figure}

In the \emph{online phase} we explore the \emph{solution manifold} for untested parameters in $\mathbb{P}$ using an approach similar to POD with Interpolation (PODI) \cite{farhat, demo2018shape}. This means that the mapping between the parameter space and the POD space is approximated by using some regression or interpolation techniques. Formally we train the regressor $\mathcal{I} : \mathbb{P} \to \mathbb{R}^\npod$ on the system of $\ntrain = |\mathcal{S}|$ vectorial equations $\mathcal{I}(\mu_i) = \widetilde{\matr{x}}_{t^*}^{\mu_i}$ for $i\in\{1, \dots, p\}$. Using $\mathcal{I}$ we are able to approximate a new reduced snapshot $\widetilde{\matr{x}}^{\mu_*}_{t^*}$ for an untested parameter $\mu_*$. We remark that in order to obtain the approximated reduced snapshot we need a multidimensional regressor; alternatively, a single regressor can be built for any components of the reduced snapshots.

At this point, we are able to recover the approximated full order vector $\matr{x}^{\mu_*}_{t^*} \approx \widehat{\matr{x}}^{\mu_*}_{t^*} = \matr{U}_n \widetilde{\matr{x}}^{\mu_*}_{t^*} \in \mathbb{R}^{\spacedim}$, which is the predicted realization of the system with the untested parameter $\mu_*$ in the future time instant $t^*$.
\RA{
We have summarizes the online procedure in Figure~\ref{fig:online}, showing the integration of the numerical tools adopted within the introduced pipeline.
}

\begin{algorithm}
    \caption{Dynamic Mode Decomposition for parametric problems | Online phase}
    \hspace*{\algorithmicindent} \textbf{Input}: The operator $\widetilde{\matr{A}}$ (or the operators $\{\widetilde{\matr{A}}^{\mu_i}\}_{i=1}^p$)\\
    \hspace*{\algorithmicindent} \textbf{Input}: the POD basis $\matr{U}_\npod$\\
    \hspace*{\algorithmicindent} \textbf{Input}: the parameter $\mu_* \in \mathbb{P}$\\
    \hspace*{\algorithmicindent} \textbf{Input}: the time instant $t^*$\\
    \label{alg:online}
    \begin{algorithmic}
        \STATE{Predict the reduced snapshots $\widetilde{\matr{x}}^{\mu_i}_{t^*} \in \mathbb{R}^\npod$ usind DMD operator(s) for $i\in\{1, \dots, p\}$}
        \STATE{Build the mapping  $\mathcal{I} : \mathbb{P} \to \mathbb{R}^\npod$ such that $\mathcal{I}(\mu_i) = \widetilde{\matr{x}}_{t^*}^{\mu_i}$ for $i\in\{1, \dots, p\}$}
        \STATE{Approximate the new reduced snapshot $\widetilde{\matr{x}}_{t^*}^{\mu_*} = \mathcal{I}(\mu_*)$}
        \STATE{Recover the full dimensional representation $\widehat{\matr{x}}_{t^*}^{\mu_*} = \matr{U}_n \widetilde{\matr{x}}^{\mu_*}_{t^*} \in \mathbb{R}^{\spacedim}$}
    \end{algorithmic}
\end{algorithm}

\section{\RA{Knowledge-based DMD operator stabilization}}
\label{meth:stab}
It may happen in some cases that DMD produces \RA{\emph{not physical}} modes \RA{in a frequency sense. The eigenvalues corresponding to the modes give us indeed precise information about their stability, which we can exploit together with our a priori knowledge about the problem at hand. Let us suppose to deal with a stable dynamical system, whose output is used to compute the DMD operator. Nonetheless the stability of the original system, DMD could compute divergent modes}, which means that the eigenvalues associated with these modes reside outside the unit circle (\emph{modulo} greater than one). This is a very problematic aspect if we intend to perform prediction on future time instants, since these eigenvalues are very likely to diverge due to Eq.~\eqref{eq:dmd_future}, making the prediction unusable.
The DMD dynamics --- i.e. the evolution in time of the modes --- is strictly correlated to the problem at hand, but it is also affected by numerical approximation, noisy data, or a limited set of input snapshots, which may induce the DMD to compute divergent modes even when applied to stable systems.
To overcome these difficulties, we propose here a simple approach to stabilize the DMD operator by looking at its eigenvalues. The main idea is to keep only the modes \RA{that show dynamics consistent with the nature of the problem, discarding all the spurious contributions. Thinking again of a stable dynamical system as an example,} the diverging and converging eigenvectors need to be neglected in order to have more robust forecasting.


We consider a threshold $\epsilon~\RA{\ll 1}$ as the unique (hyper)parameters of such a stabilization. First of all, we discard all the DMD modes that have the modulo of the corresponding eigenvalue greater than $1 + \epsilon$, which we remark be the divergent ones.
We perform the same procedure also to convergent modes --- which tend to have a nil contribution advancing in time ---, removing the modes that have the corresponding eigenvalue inside the unit circle.
\begin{figure}
    \centering
    \includegraphics[width=.8\textwidth]{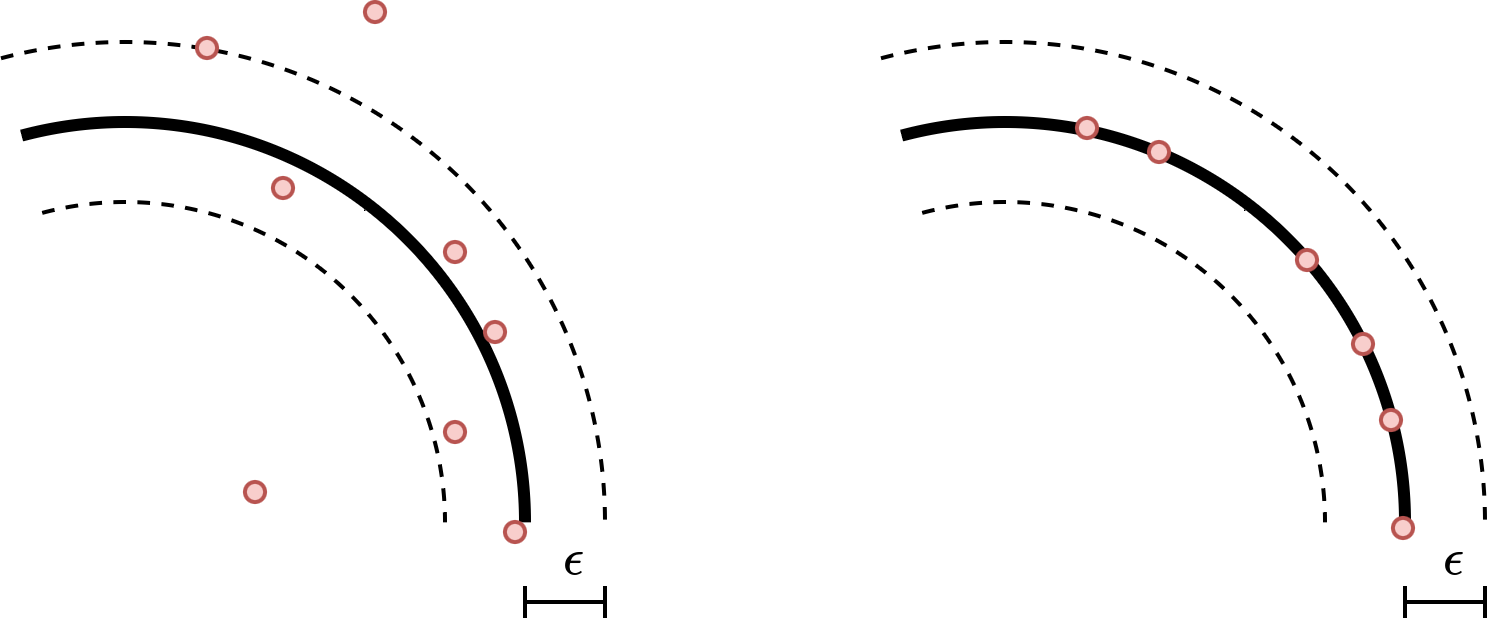}
    \caption{\RA{Graphical example of stabilization: on left, the eigenvalues (red) of the DMD operator are sketched together with the unitary circle (solid black line) in the complex plane. On right, we show the same entities after the stabilization: the eigenvalues distant from the unit circle have been discarded, whereas the closest are scaled in order to have the modulo equal to one.}}
    \label{fig:stab_dmd}
\end{figure}
Finally, the eigenvalues of the remaining modes are projected onto the unity circle (i.e. the eigenvalue is normalized), in order to stabilize the operator.
Formally, we manipulate the eigenvalues such that:
\begin{equation}
\lambda_i^\text{stable} =
\begin{cases}
0  & \text{if}\,\mid \lambda_i - 1 \mid > \epsilon\\
\frac{\lambda_i}{\mid \lambda_i \mid}  & \text{otherwise}\\
\end{cases}\quad\quad\quad
\text{for}\,i=1, \dots, r.
\end{equation}
Defining $\mathbf \Lambda = \diag{(\lambda_1^\text{stable}, \dots, \lambda_r^\text{stable})}$ is sufficient to obtain the stabilized prediction since
\RA{
the contribution of unstable modes is neglected by applying Eq.~\eqref{eq:dmd_future}. It must be said that this discarded contribution we apply to the DMD operator may cause a higher error in the reconstruction of the training snapshots, but it aims to limit such an error in unseen time instants, resulting in an ideal procedure to apply for future forecasting.
Finally, we conclude by specifying the generality of this approach: during one of the experiments collected in this work, we applied such processing to stabilize the DMD operator(s), since we knew the periodically stable nature of the studied problem (\ref{sec:cylinder}). The same approach could be applied by exploiting all the knowledge we have: in the case of a dumped system, only converging (eigenvalues modulo less than one) modes should be considered, whereas if dominant frequencies should be present, only modes with those frequencies (or their multiples) could be taken into account to improve the future accuracy,
}

\section{Numerical results}
\label{sec:results}
In this section we present the numerical results obtained using the two approaches discussed in \cref{meth:full,meth:partial} on three problems in different topics.

\begin{figure}[t!]
    \begin{subfigure}[b]{0.30\textwidth}
        \centering
        \includegraphics[width=\textwidth]{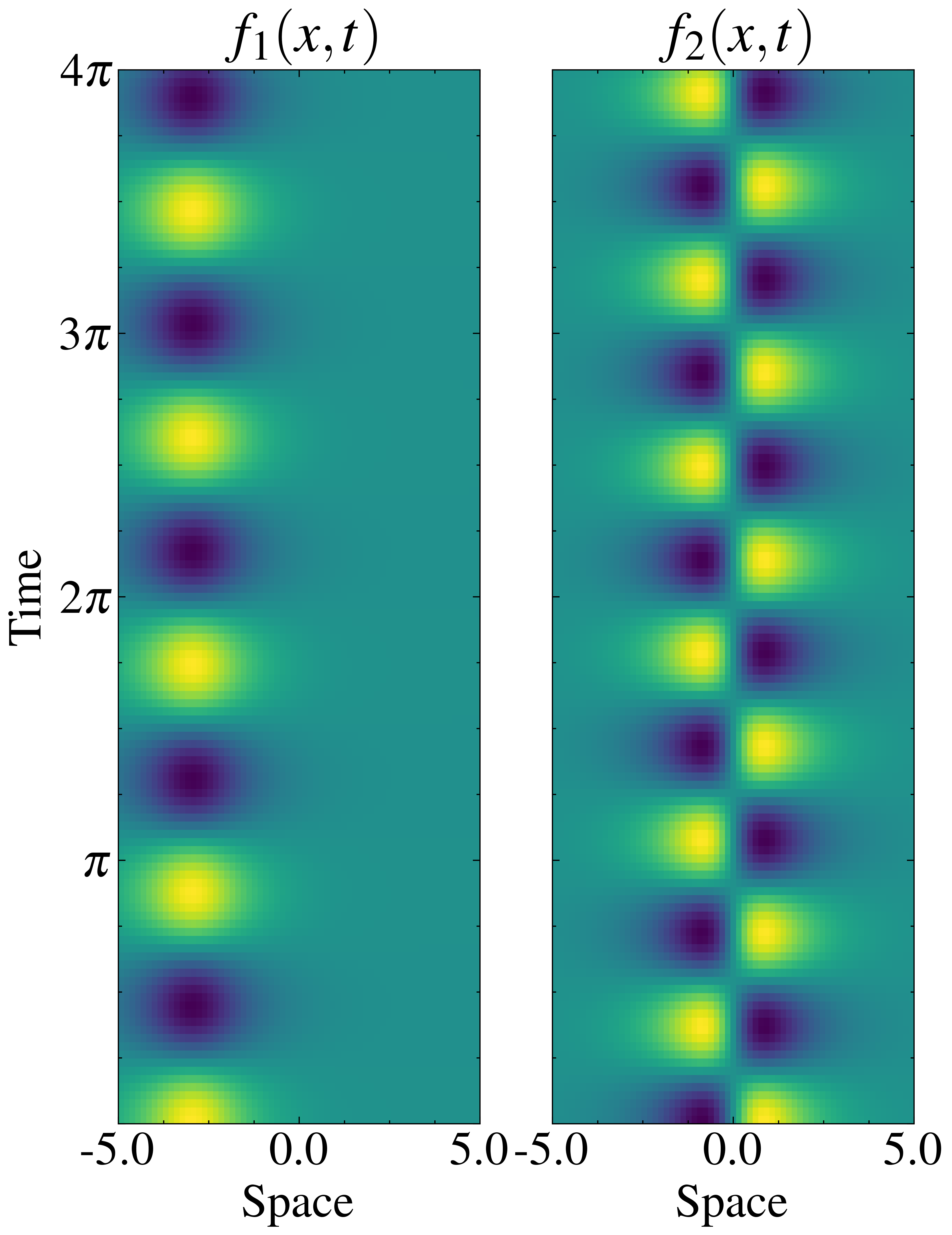}
    \end{subfigure}
    \hfill
    \begin{subfigure}[b]{0.68\textwidth}
        \centering
        \includegraphics[width=\textwidth]{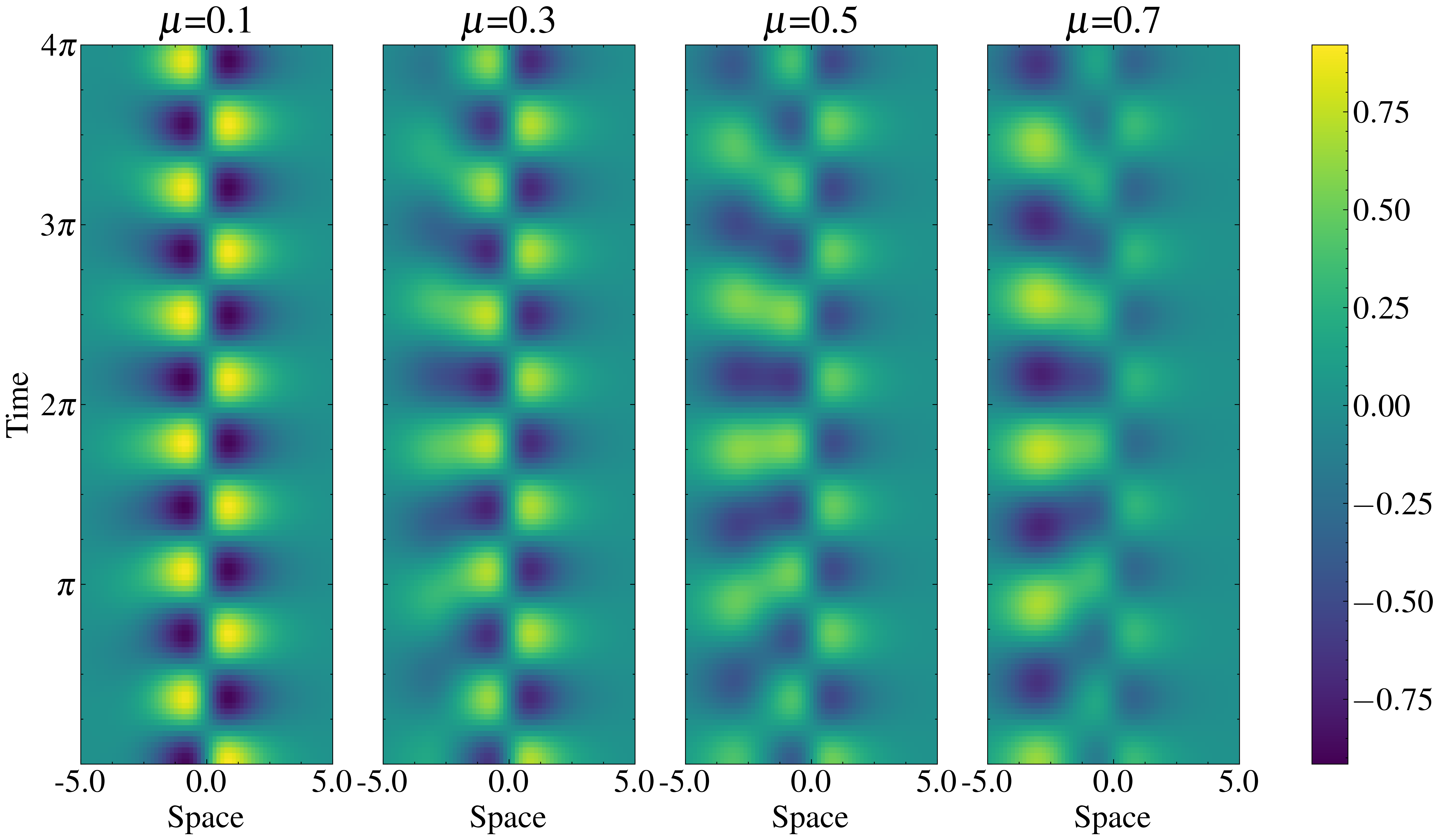}
    \end{subfigure}
    \caption{On the left, the two simple functions which we are using for the composition. We consider 10 values of $\mu$ to train our algorithm, the realizations corresponding to $\mu \in \{0.1,0.3,0.5,0.7\}$ are shown on the right.}
    \label{toy:system_showcase}
\end{figure}

\begin{figure}[b!]
    \centering
    \begin{adjustbox}{minipage=\linewidth,scale=0.8}
        \begin{subfigure}[b]{0.78\textwidth}
            \includegraphics{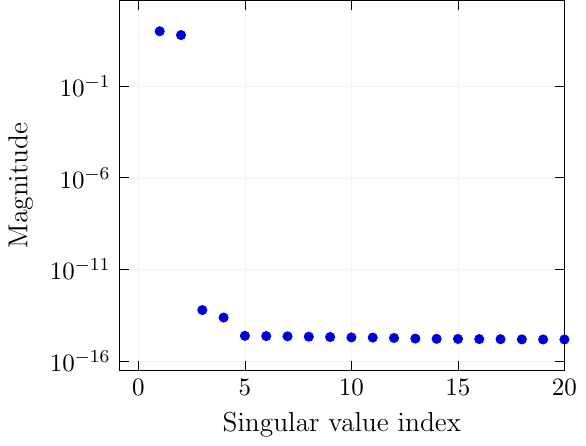}
            \caption{First 20 singular values extracted from $\mathcal{X}_1$.}
            \label{toy:pod_sv}
        \end{subfigure}
        \hfill
        \begin{subfigure}[b]{0.2\textwidth}
            \centering
            \includegraphics[angle=90,origin=c, width=0.985\textwidth]{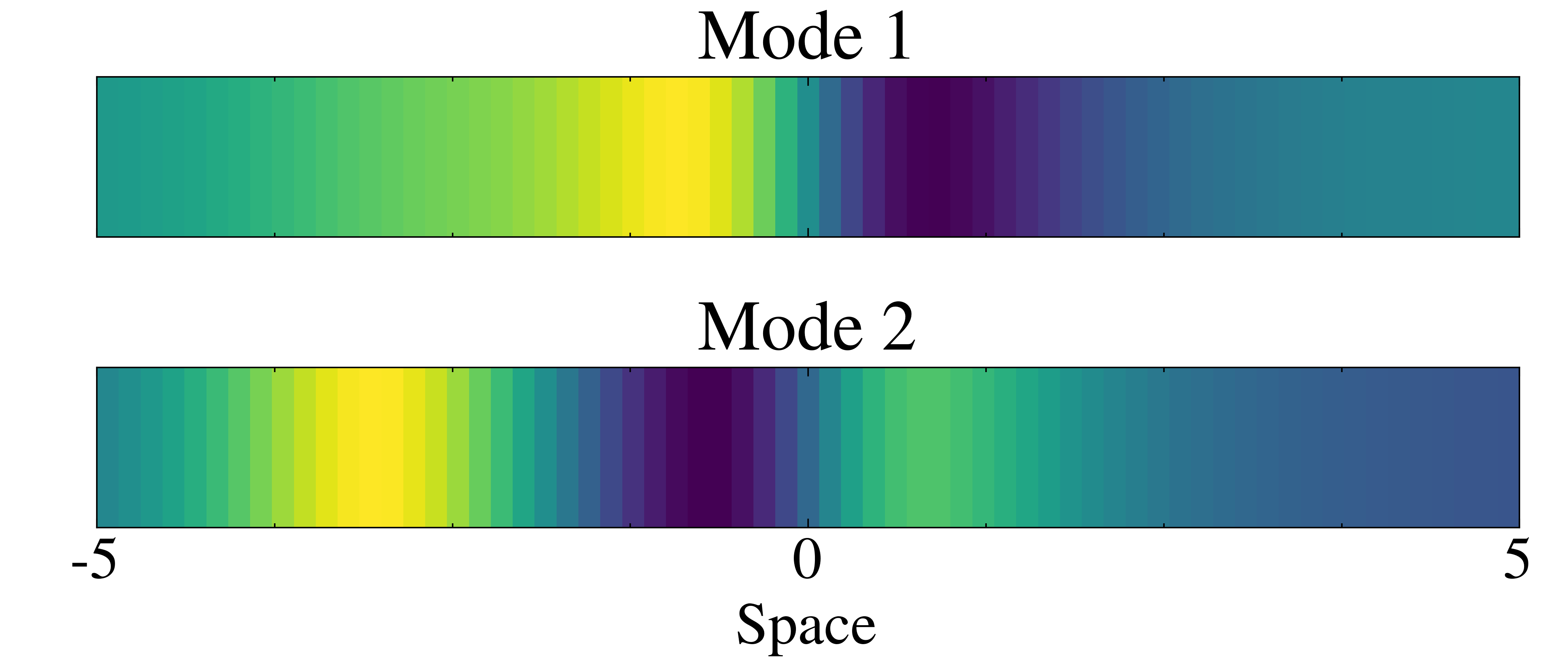}\hfill
            \caption{First two POD modes extracted from $\mathcal{X}_1$.}
            \label{toy:pod_modes}
        \end{subfigure}
    \end{adjustbox}
    \caption{Looking at singular values of the system, on the left, we deduce that the most contribution is given by the first two POD modes. We drew those modes on the right.}
\end{figure}

\subsection{Numerical study 1: System of algebraic equations with linear parameter dependency}
\label{sec:toy}

For our first experiment we consider the simple system depicted in \cref{toy:system_showcase}, which is the parametrized composition of two time dependent functions:
\begin{gather}\label{res:first_nonsampled}
    f(x,t,\mu) = \mu f_1(x,t) + (1-\mu) f_2(x,t), \quad \begin{cases}
        f_1(x, t) = \sech (x+3)\exp (i2.3t)\\
        f_2(x, t) = 2\sech (x)\tanh(x) \exp(i2.8t)
    \end{cases}
\end{gather}
We consider the time instants set $\mathcal{T}$ to be an equispaced sampling of the time frame $[0, 4\pi]$ (we retain $\timedim=129$ samples); the spatial domain $\Omega$ is $[-5, 5]$ (we take $\spacedim=1000$ equispaced samples from this interval). Therefore our output function, which is a discrete version of that shown in \cref{res:first_nonsampled}, has the form $\matr{f}_k^\mu \in \mathbb{R}^{1000}, k \in \mathcal{T}$: each component of the vector corresponds to one of the 1000 samples of the space domain $\Omega$. From the formulation of the problem it is clear that the dependency of the system on the parameter $\mu$ is linear. We consider the training parameters set $\mathcal{S}$ to be an equispaced sampling of the set $[0,0.9]$ (we retain $\ntrain=10$ samples from this interval).

\begin{figure}[t!]
    \begin{subfigure}[b]{0.52\textwidth}
        \includegraphics{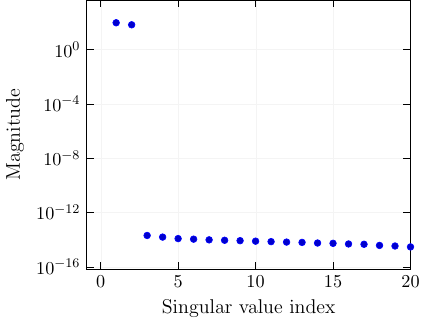}
        \caption{First 20 singular values extracted during the application of DMD to $\matr{\mathcal{X}}_2$.}
        \label{toy:dmd_sv}
    \end{subfigure}
    \hfill
    \begin{subfigure}[b]{0.39\textwidth}
        \centering
        \includegraphics[width=0.98\textwidth]{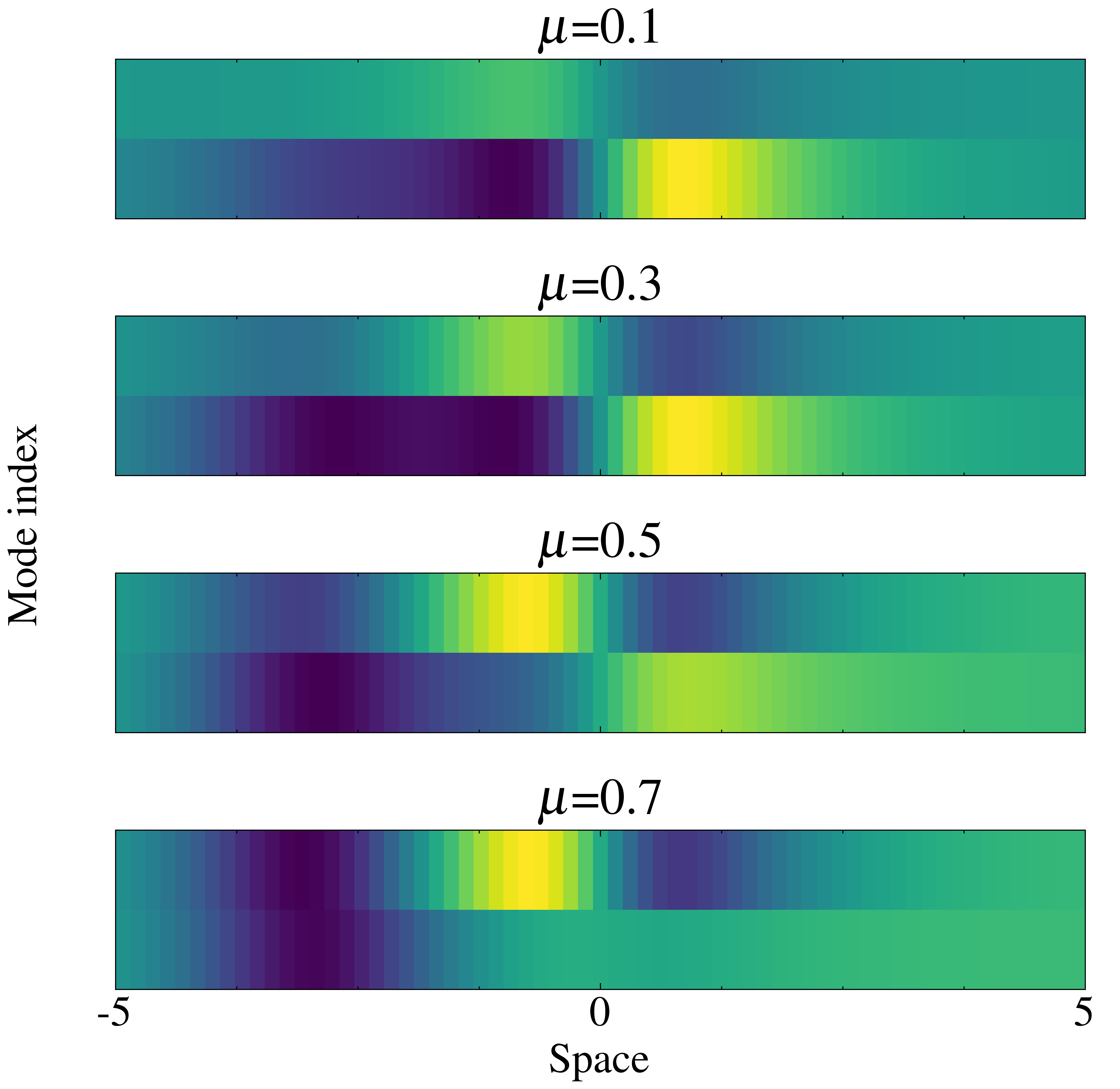}
        \caption{Appearence of the two DMD modes extracted from $\matr{\mathcal{X}}_2$, divided into sub-components.}
        \label{toy:dmd_modes}
    \end{subfigure}
    \begin{subfigure}[b]{\textwidth}
        \caption{Dynamics extracted with DMD applied on $\matr{\mathcal{X}}_2$.}
        \includegraphics{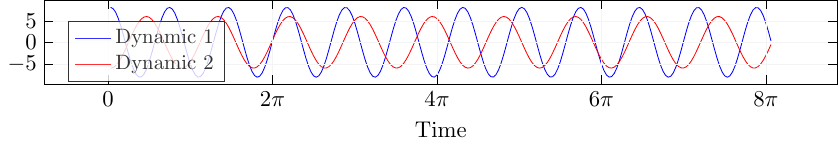}
        \label{toy:dmd_dynamics}
    \end{subfigure}
    \caption{Results of the application of DMD to the matrix $\matr{\mathcal{X}}_2$: singular values (\cref*{toy:dmd_sv}), DMD modes (\cref*{toy:dmd_modes}), dynamics of the extracted modes (\cref*{toy:dmd_dynamics}).}
\end{figure}

\begin{figure}[t!]
    \centering
    \begin{subfigure}[t]{0.49\textwidth}
        \includegraphics{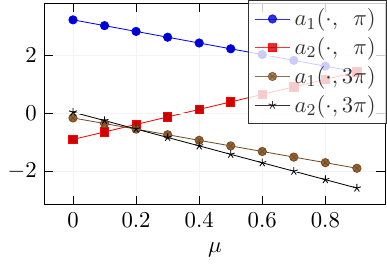}
        \caption{First and second POD coefficients for \emph{tested} parameters at several time instants.}
        \label{toy:pod_coeffs}
    \end{subfigure}
    \hfill
    \begin{subfigure}[t]{0.49\textwidth}
        \includegraphics{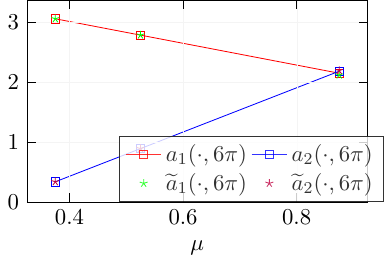}
        \caption{Interpolated/predicted and expected POD coefficients at $t=6\pi$.}
        \label{toy:expected_pod_coeffs}
    \end{subfigure}
    \caption{On the left we show the arrangement of POD coefficients for $t = \pi,3\pi < \timedim$ as $\mu$ varies. On the right, the results of the prediction and interpolation of five untest ed parameters for $t = 6\pi > \timedim$.}
\end{figure}

\begin{figure}[b!]
\begin{subfigure}[t]{0.49\textwidth}
    \centering
    \includegraphics[width=\textwidth]{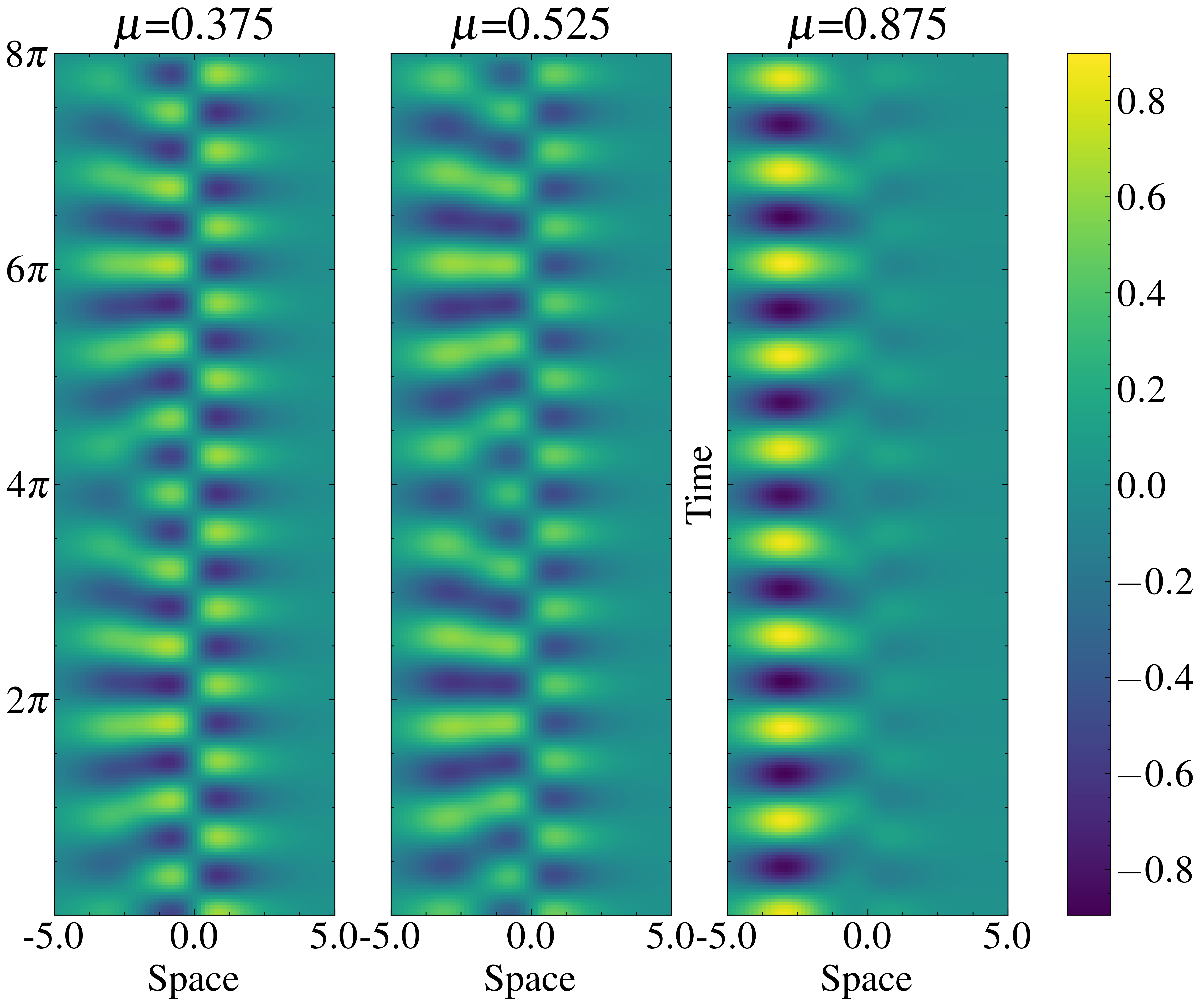}
    \caption{Reconstruction of the system for untested parameters in future time instants.}
    \label{toy:pdmd_result}
\end{subfigure}
\hfill
\begin{subfigure}[t]{0.49\textwidth}
    \centering
    \includegraphics[width=\textwidth]{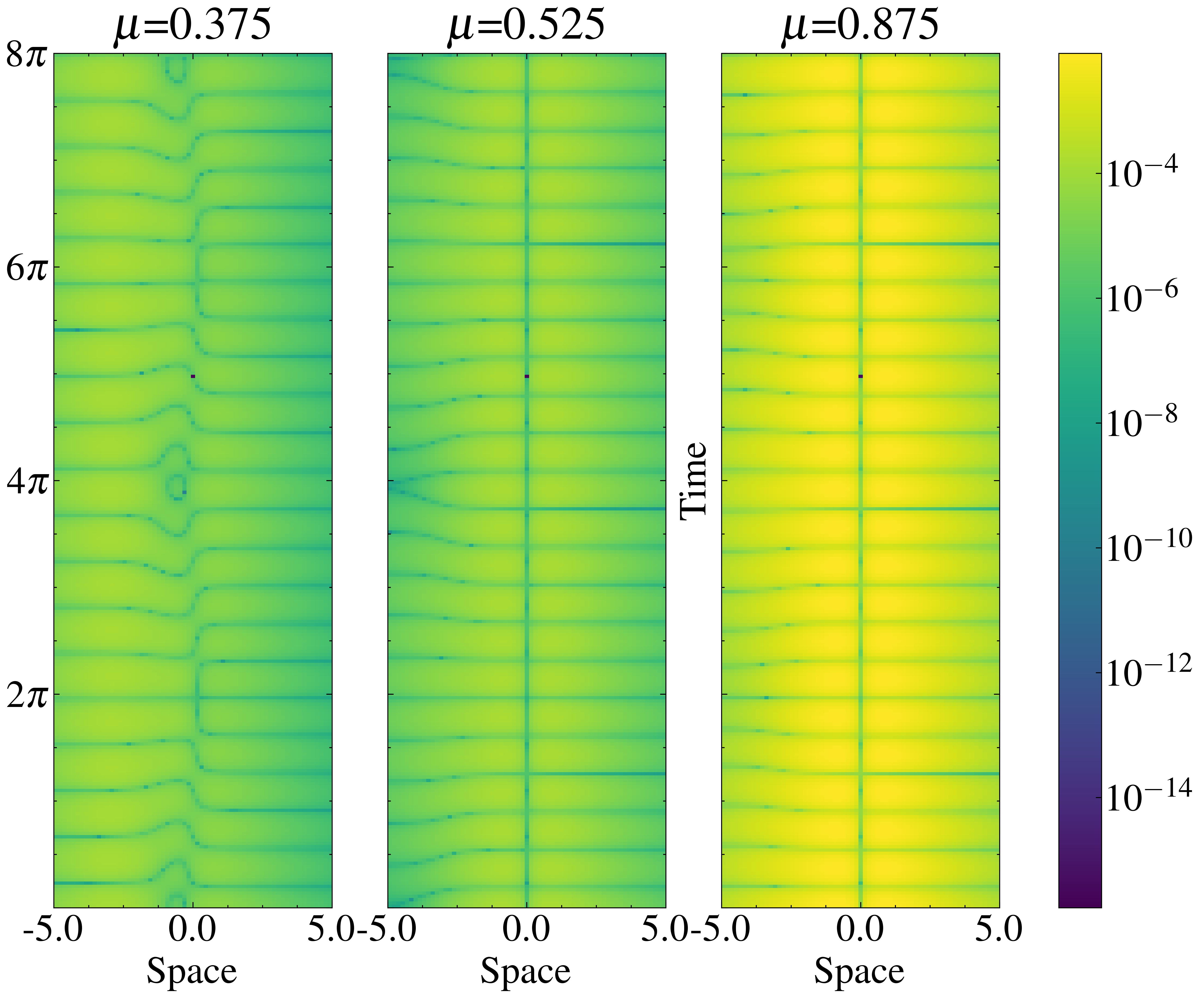}
    \caption{Absolute error of the reconstructed system.}
    \label{toy:pdmd_error}
\end{subfigure}
\caption{On the left we plot the result of the algorithm discussed in \cref{meth:full}; on the right we show the absolute error between the expected and reconstructed system.}
\label{toy:results}
\end{figure}

Considering the linear dependency of the system on the parameter $\mu$, we apply the monolithic method discussed in \cref{meth:full}. From \cref{toy:pod_sv} we infer that the optimal value for $\npod$ (the number of POD modes retained for each realization of the system) is $\npod = 2$. The optimality of this value is also confirmed by the fact that the corresponding POD modes (\cref{toy:pod_modes}) are physically coherent with the system depicted in \cref{toy:system_showcase}. After the application of DMD (we take only the first two modes employing the singular values criteria, \cref{toy:dmd_sv}) we obtain the matrix $\matr{\Phi} \in \mathbb{R}^{\npod \ntrain \times 2} = \mathbb{R}^{20 \times 2}$ of the two DMD modes extracted from $\matr{\mathcal{X}}_2$ defined like in \cref{eq:pod_coefficients_matrix}. In order to visualize the results of the analysis we split $\matr{\Phi}$ in $\ntrain$ sub-components $\matr{\Phi}_i \in \mathbb{R}^{\npod \times 2} = \mathbb{R}^{2 \times 2}$: each column contain the part of the DMD mode relative to a particular parameter. After that we use the matrix of POD modes to recover the full-dimension appearance of each DMD mode (\cref{toy:dmd_modes}). The dynamics of the two DMD modes are then visualized in \cref{toy:dmd_dynamics}, as expected they are periodic functions.

We are now able to start the \emph{online phase}: using DMD modes and dynamics we approximate POD coefficients corresponding to future time instants ($t > \timedim$) for tested values of the parameter. We would like to double the amplitude of the time frame, namely we intend to approximate the behavior of the system until $t = 8\pi$ (we keep the same time difference between consecutive snapshots). We also consider three untested parameters: $\mu_x \in \{0.375, 0.525, 0.875\}$. After the examination of \cref{toy:pod_coeffs} we can assume that POD coefficients are arranged roughly along a straight line for all time instants. This situation is ideal for the interpolation, assuming that DMD preserves this kind of arrangement. In \cref{toy:expected_pod_coeffs} we plot the results of the interpolation, versus the expected POD coefficients, which we obtained projecting the matrix $\matr{X}^{\mu_x}$ onto the reduced POD subspace spanned by the columns of the matrix $\matr{U}_\npod$ of POD modes, where $\mu_x$ is any of the untested parameters which we mentioned above.

In \cref{toy:results} we depict the system predicted/interpolated using our approach, for the three untested parameters $\mu_x \in \{0.375, 0.525, 0.875\}$, in the doubled time frame $[0,8\pi]$. The non-null error is most likely given by the fact that we considered a small number of POD modes with respect to the high dimensionality of the original system.


\subsection{Numerical study 2: An heat conductivity problem}
\label{sec:heat}
We consider now a parametric nonlinear unsteady heat problem introduced in~\cite{grepl2007efficient,hoang2021projection}. We want to compute the unknown $u(\mathbf x, t, \mupar)$ with $\mathbf x = (x_1, x_2) \in \Omega = [0, 1] \times [0, 1]$, $t \in [0, 2]$ and $\mupar = (\mu_1, \mu_2) \in \mathbb{P} = [0.01, 10]^2$ such that
\begin{equation}
    \begin{cases}
    \frac{\partial u}{\partial t} - \nabla^2 u = 100\sin{(2\pi x_1)}\sin{(2\pi x_2)}\sin{(2\pi t) - \frac{\mu_1}{\mu_2}(e^{\mu_2u} -1)} & \text{in}\,\Omega,\\
    u = 0 & \text{on}\,\partial \Omega.
    \end{cases}
\end{equation}

The discrete solution of such a problem is obtained by employing the finite element framework. The 2D spatial domain is divided in $800$ triangular cells with P1 element.

\begin{figure}[t!]
    \begin{subfigure}[b]{0.477\textwidth}
        \includegraphics{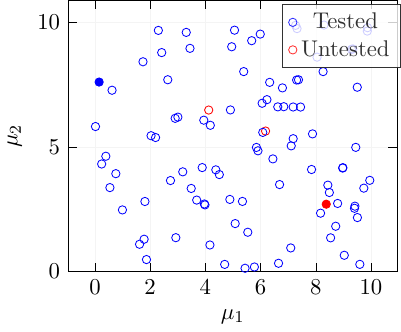}
        \caption{Distribution of tested and untested parameters for the second test case.}
        \label{heat:params}
    \end{subfigure}
    \hfill
    \begin{subfigure}[b]{0.47\textwidth}
        \includegraphics[width=\textwidth]{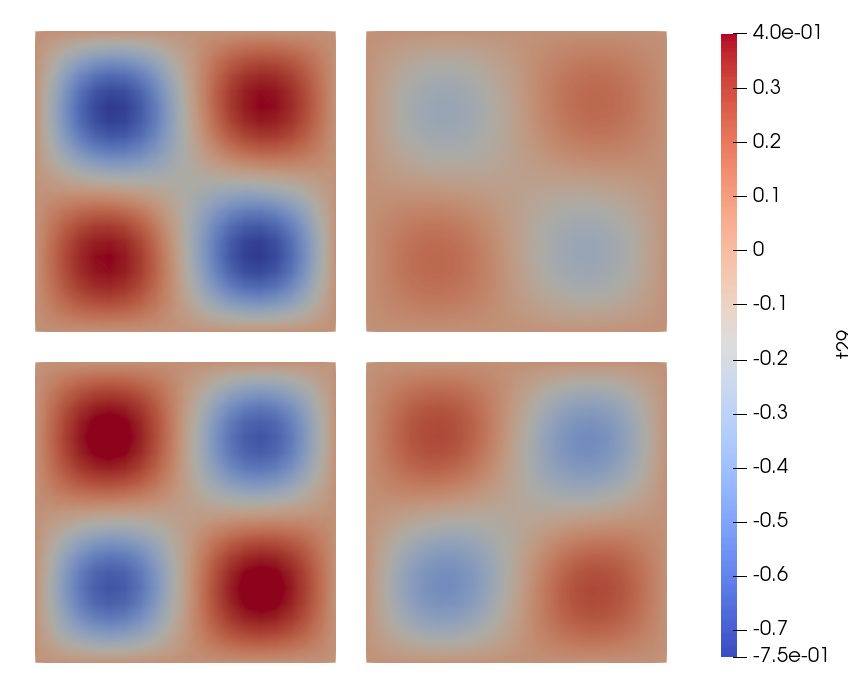}
        \caption{The system for the $\mu$ in solid blue in the left figure, sampled at four time instants.}
        \label{heat:system}
    \end{subfigure}
    \caption{\RB{On the left we plot the distribution of parameters we considered for the system, and on the right we depict the realization of the system represented by the solid blue point in four time instants.}}
\end{figure}

\begin{figure}[t!]
    \begin{subfigure}[b]{0.49\textwidth}
        \includegraphics{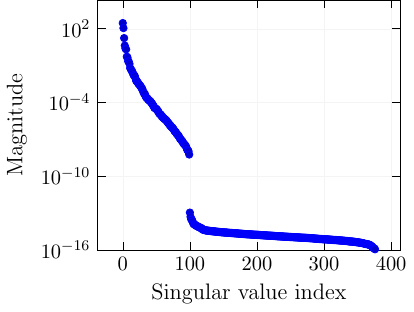}
        \caption{Singular values extracted from $\mathcal{X}_1$.}
        \label{heat:pod_sv}
    \end{subfigure}
    \hfill
    \begin{subfigure}[b]{0.49\textwidth}
        \includegraphics[width=\textwidth]{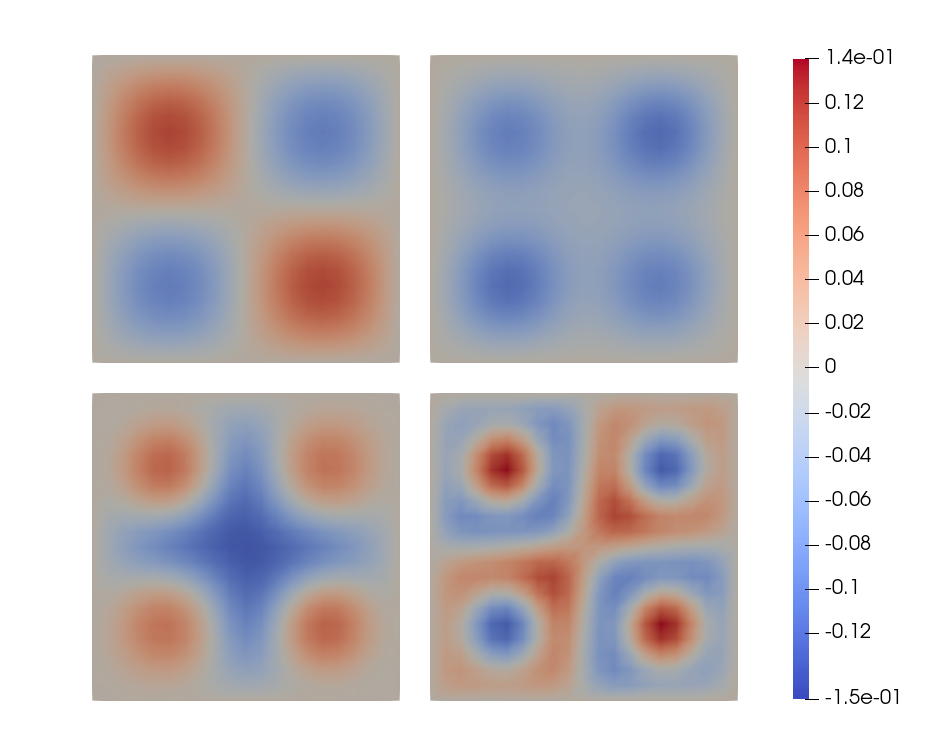}
        \caption{First four POD modes extracted from $\mathcal{X}_1$.}
        \label{heat:pod_modes}
    \end{subfigure}
    \caption{Unlike the system shown in \cref{sec:toy}, it is clear that we would need to consider an higher number of POD modes in this case in order to guarantee an high fidelity reconstruction.}
\end{figure}

We consider 95 instances of the system. The distribution of the parameters (tested and untested) is depicted in \cref{heat:params}, alongside the physical appearance of the system at several (known) time instants.
For all the parametric samples, we collect $86$ snapshots equispaced in the time window $[0, 85]$.
We retained 30 POD modes for the reduced space, and did not perform an additional reduction during the DMD-step. We attempt to predict/interpolate the system for three more \emph{untested} parameters, depicted in red in \cref{heat:params}, which we collect into the set $\mathcal{Q}$. We also consider an extended time frame $86 \leq t \leq 100$.

In \cref{heat:pod_sv} we observe that, unlike the system considered for the first experiment, POD applied on $\matr{\mathcal{X}}_1$ needs a considerable number of POD modes in order to capture all the significative characteristic structures of the parametric system. The first hard step occurs around the 100-th singular value, however we preferred to choose a lower threshold in order to keep the interpolation error under control. We can still see that the first modes are the most significative ones, even when POD is applied simultaneously on more realizations of the same parametric system.


\begin{figure}[b!]
    \begin{subfigure}[t]{0.49\textwidth}
        \includegraphics{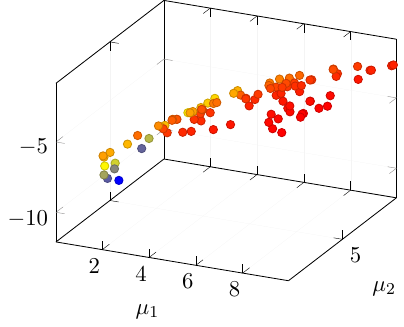}
        \caption{First POD coefficient at $t=10$.}
        \label{heat:pod_coeffs_1_t10}
    \end{subfigure}
    \hfill
    \begin{subfigure}[t]{0.49\textwidth}
        \includegraphics{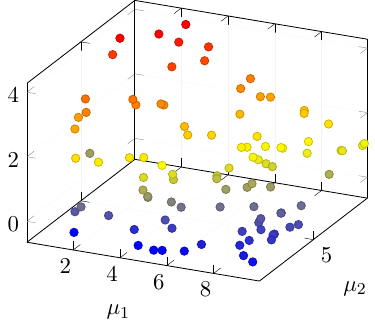}
        \caption{Second POD coefficient at $t=10$.}
        \label{heat:pod_coeffs_2_t10}
    \end{subfigure}
    \caption{The first POD coefficients (\cref{heat:pod_coeffs_1_t10}) for each value of $\mu$ are arranged roughly on a flat surface. However, the second POD coefficients (\cref{heat:pod_coeffs_2_t10}) are much more scattered, but we could attempt the interpolation using a linear interpolator, or the algorithm K-NearestNeighbors.}
    \label{heat:pod_coeffs}
\end{figure}

\noindent In \cref{heat:pod_coeffs} we see that the distribution of POD modal coefficients (starting from the second one) is less regular than the distribution in the first experiment.

\begin{figure}[t!]
    \centering
    \includegraphics[width=\textwidth]{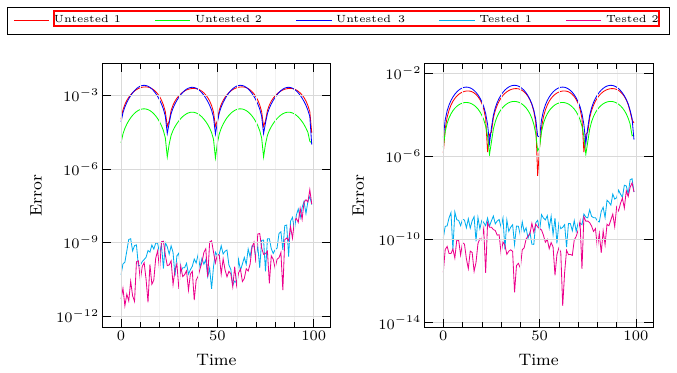}
    \caption{Evolution over time of the absolute error on the first and second POD coefficient for several \emph{tested} and \emph{untested} values of the parameter.}
    \label{heat:pod_error}
\end{figure}

\begin{figure}[b!]
    \begin{subfigure}[t]{0.43\textwidth}
        \includegraphics[width=\textwidth]{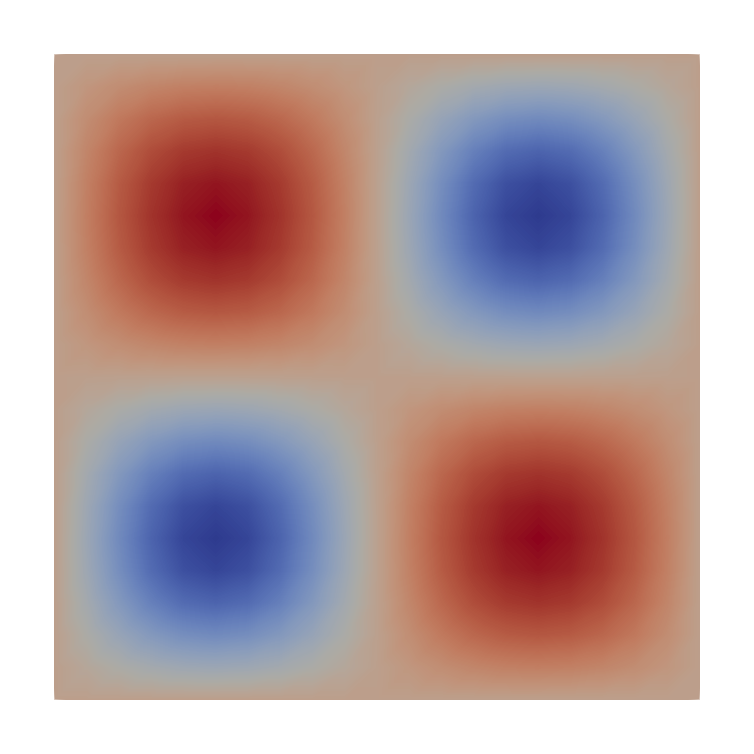}
        \caption{Original system.}
        \label{heat:original_untested_system}
    \end{subfigure}
    \hfill
    \begin{subfigure}[t]{0.545\textwidth}
        \includegraphics[width=\textwidth]{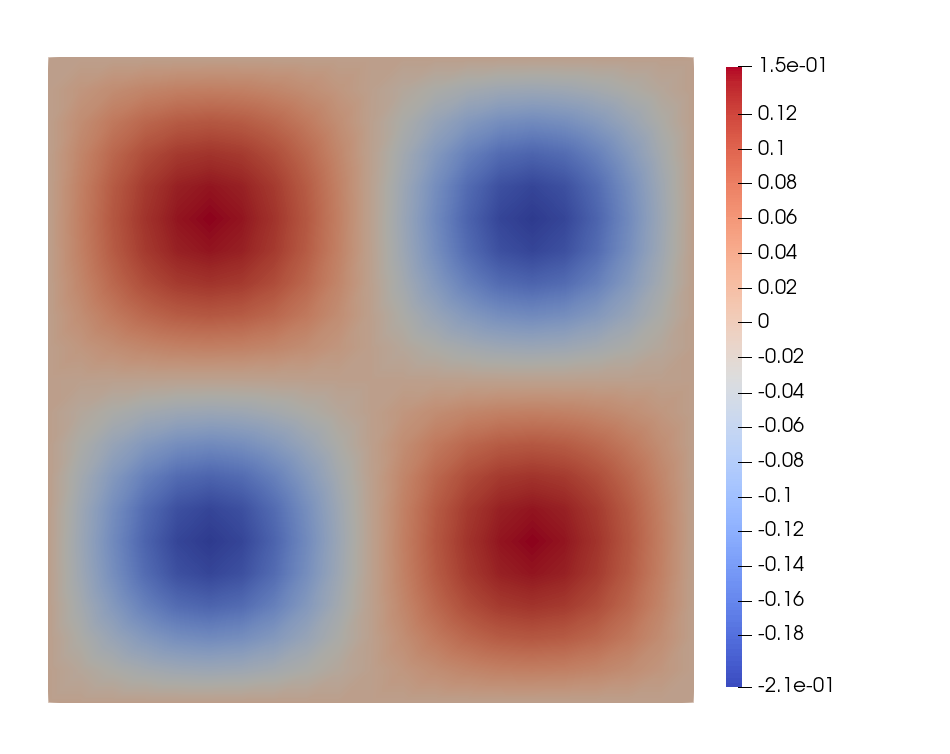}
        \caption{Output returned by pDMD.}
        \label{heat:pdmd_reconstruction}
    \end{subfigure}
    \label{heat:pdmd_results}
    \caption{Graphical result for an untested parameter in the heat problem, at $t=91$.}
\end{figure}
\begin{figure}[t!]
    \begin{subfigure}[t]{0.49\textwidth}
        \includegraphics[width=\textwidth]{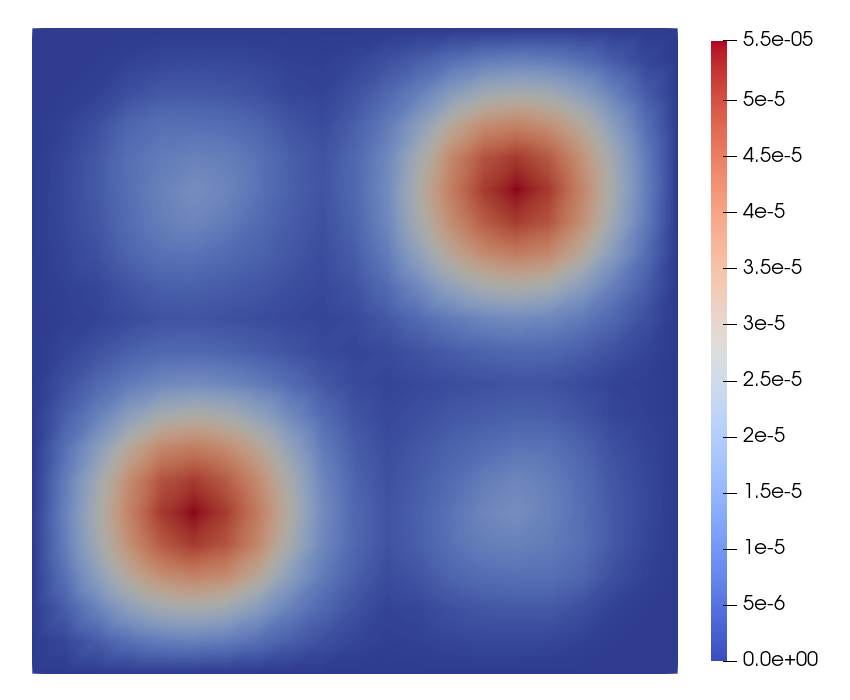}
        \caption{Absolute error.}
        \label{heat:pdmd_absolute_error}
    \end{subfigure}
    \hfill
    \begin{subfigure}[t]{0.49\textwidth}
        \includegraphics[width=\textwidth]{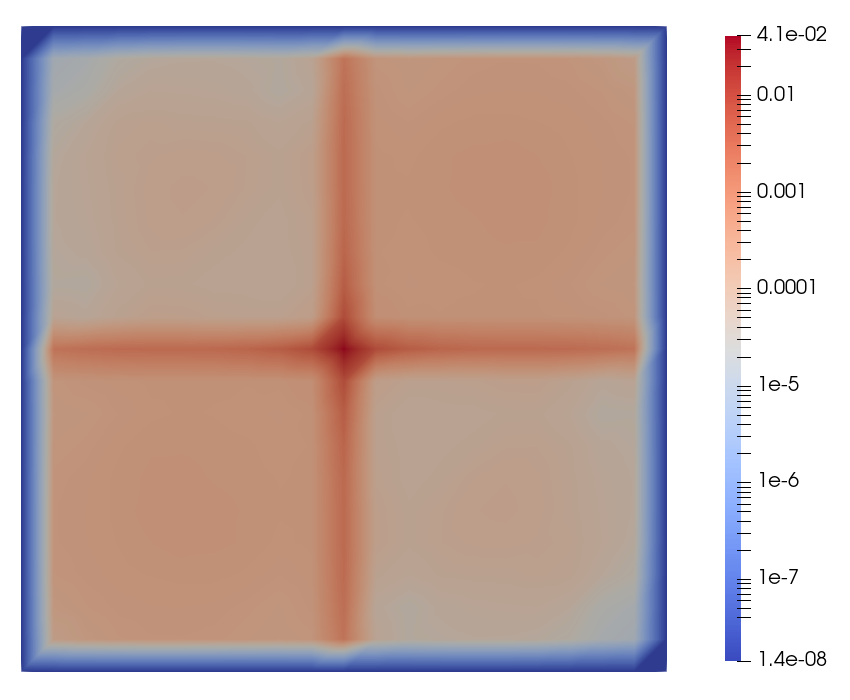}
        \caption{Relative error.}
        \label{heat:pdmd_relative_error}
    \end{subfigure}
    \caption{Error plot for the untested parameter (same as \cref{heat:pdmd_results}}.
\end{figure}

\noindent In \cref{heat:pod_error} we plot the result of the interpolation (we used a linear interpolator) for two \emph{tested} and three \emph{untested} parameters. The error on POD coefficients for tested parameters are visualized to check in which step we introduce the biggest error. In fact we can see that in the unknown time frame (starting from $t = 85$) the error on tested POD coefficients starts to increase slightly, while the error on untested POD coefficients remains stable. This means that the error introduced with the prediction provided by DMD is negligible with respect to the error introduced by interpolators.

Finally, we plot the prediction for one of the \emph{untested} parameters and unknown time instants computed using our approach, alongside the truth solution, the relative and absolute error. The red cross at the center of \cref{heat:pdmd_relative_error} is given by the fact that that region of $\Omega$ maintains values around zero for the majority of time instants, therefore the relative error appears high with respect to the neighboroods. The same ``problem'' does not appear in \cref{heat:pdmd_absolute_error}, since no division is performed.

For this problem we consider a set of experiments aimed at checking the sensitivity of the approach on the number of tested parameters and on the size of the time frame used during the training of DMD. We also take into account several interpolators with different features for the online phase, in order to assess which one could vcally provide better performances. In order to visualize such measurements we are going to compute a scalar function $e_\mathcal{I}(\mathcal{S}, \mathcal{T}, t, \mu)$ defined as the relative error at time $t$ between the reconstructions of the system in the untested parameter $\mu \in \mathcal{Q}$ and the high-fidelity validations.
Since we consider more than one $\mu$ for this measurement (we recall that $\mathcal{Q}$ is a set of untested parameters which we presented above), we are going to consider the mean relative error. Formally we define it as:
\begin{gather}
    e_\mathcal{I}(\mathcal{S}, \mathcal{T}, t) = \frac{1}{|\mathcal{Q}|}\sum_{\mu \in \mathcal{Q}} \frac{\Vert \hat{\matr{x}}_{t}^\mu[\mathcal{S}, \mathcal{T}] - \matr{x}_t^\mu \Vert_2}{\Vert \matr{x}_t^\mu\Vert_2}
\end{gather}
where $\hat{\matr{x}}_{t}^\mu[\mathcal{S}, \mathcal{T}]$ refers to the snapshots approximation of the algorithm trained over $\mathcal{S}$ parametric samples and $\mathcal{T}$ time instants.
 We would like to investigate the accuracy of the prediction in the experiments, therefore we are going to evaluate the value of $e_\mathcal{I}(\mathcal{S}, \mathcal{T}, t)$ for several combinations of the parameters taken into account, in order to assess the sensitivity of the proposed method on the different factors which affect the quality of the reconstruction.

Since we need to interpolate values on a 2-dimensional surface, during the experiment we are going to consider the following regressors, for which we provide a brief explanation:
\begin{itemize}
    \item \emph{Linear}: divides the surface in 2-dimensional simplices, and interpolates linearly on each sub-surface \cite{hoang2021projection};
    \item \emph{Cubic}: return the value determined from a piecewise cubic polynomial surface \cite{hoang2021projection};
    \item \emph{Nearest neighbor}: Return the nearest value in the dataset \cite{nearest};
    \item \emph{GPR}: Interpolates using a Gaussian Process \cite{gaussian_process}.
\end{itemize}

In the first experiment we consider a random initial training set $\mathcal{S}_0$ which contains 20 parameters. We add random parameters one by one in such a way that $\mathcal{S}_{k+1} = \mathcal{S}_k \cup \texttt{rand}(\mathbb{P} - \mathcal{S}_k)$, where the function $\texttt{rand}$ extracts a random member from the set passed as argument, and $\mathbb{P}$ is the set which contains all the allowed parameters for the parametric system taken into account. This procedure is iterated $r$ times, until the size of the training set (which we refer to as $\ntrain=|\mathcal{S}_{r}|$) is equal to $95$. In this experiment the size of the time frame used for the training phase is fixed ($\mathcal{T} = \{t \mid 1 \leq t \leq 85\}$). We also consider five distinct untested parameters $\widehat{\mathcal{S}}$ extracted randomly from the set $\mathbb{P} - \mathcal{S}_r$. The distribution of parameters used for the training phase and those used to validate the approach is plotted in \cref{heat:params}.

\begin{figure}[t!]
    \centering
    \includegraphics[width=\textwidth]{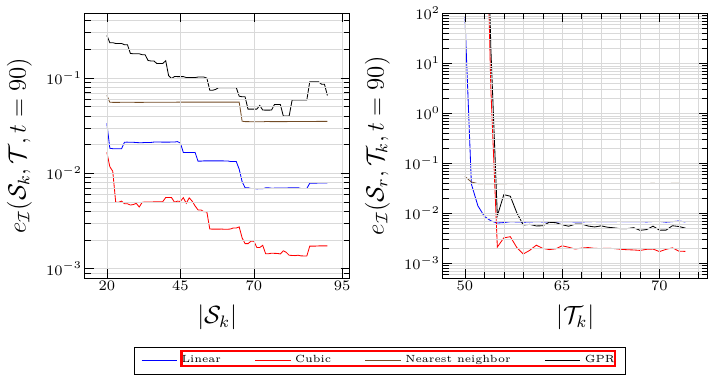}
    \caption{\RB{Sensitivity on the number of tested parameters ($\mathcal{S}_k$) and on the number of tested time instants ($\mathcal{T}_k$)}}.
    \label{sensitivity}
\end{figure}

Using the approach presented in \cref{meth:partial}, for each training set $\mathcal{S}_k$ we interpolate the system in the untested parameters $\widehat{U}$ and predict the system in the time frame $86 \leq t \leq 100$. We also evaluate the function $e_\mathcal{I}(\mathcal{S}_k, \mathcal{T}, t)$ for each possible combination of its arguments (keeping $\mathcal{T}$ fixed, and $k \in \{0, \dots, r\}$). The value $e_\mathcal{I}(\mathcal{S}_k, \mathcal{T}, t=90)$ is shown in the left part of \cref{sensitivity}: we depict the dependency of the mean relative error on the cardinality of the set $\mathcal{S}_k$ of training time instants considered during the offline phase. It is clear that the error decreases in general when the number of tested parameter increases as expected; however there is a little bit of \emph{overfitting} around $\ntrain = 90$ for all the interpolators.

In the second experiment the set of tested parameters is fixed ($\ntrain=95$, we consider the biggest training set $\mathcal{S}_r$ from the last experiment). We consider an initial time frame $\mathcal{T}_0$ which spans from 1 to 50; at each step we add one time instant to the training time frame in such a way that $\mathcal{T}_{k+1} = \min(\mathbb{N} - \mathcal{T}_k)$. We iterate this process $q$ times, until $85 \in \mathcal{T}_r$.

Then we execute the same operations performed in the last experiment: we interpolate and predict the system in the time frame $86 \leq t \leq 100$ for five untested parameter (the same of the last experiment), and we compute the values of the function $e_\mathcal{I}(\mathcal{S}_r, \mathcal{T}_k, t)$ for all the possible combinations of its parameters (this time we keep $\mathcal{S}_r$ fixed and let $k \in \{0, \dots, q\}$).

The results of the experiment are plotted in the right part of \cref{sensitivity} for $t=90$, where we depict the dependency of the error function on the cardinality of the set $\mathcal{T}_k$ of training time instants considered during the offline phase. After a critical phase, which is most likely due to the internal working of DMD (we remark that the error taken into account is relative), the error decades as the number of time instants increases, which is the expected behavior; finally (time frame $0 \leq t \leq 60$, approximately) the error becomes almost stable.


\subsection{Numerical study 3: Navier-Stokes's cylinder}\label{sec:cylinder}
The third numerical experiment is the simulation of a parametric Navier-Stokes flow passing around a circular cylinder. The equations that describes the problem are:
\begin{figure}[t!]
    \centering
    \includegraphics{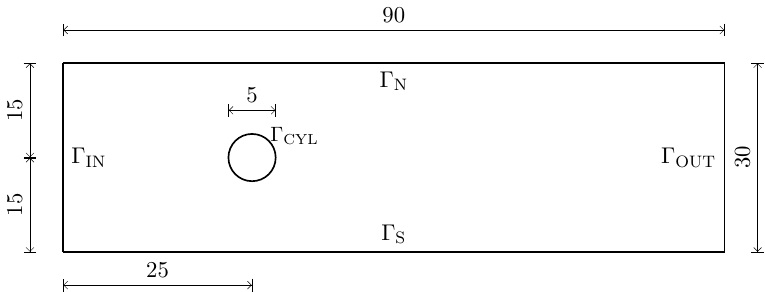}
    \caption{The domain $\Omega$ for the parametric Stokes flow simulation.}\label{fig:stokes_domain}
\end{figure}
\begin{equation}
\begin{cases}
\frac{\partial u}{\partial t} + (u \cdot \nabla)u
+\nabla p - \nu \Delta u=f & \text{ on } \Omega \times [0, T],\\
\nabla \cdot u=0 & \text{ on } \Omega \times [0, T],\\
u_x = \mu & \text{ in } \Gamma_\text{IN} \times [0, T],\\
u = 0 & \text{ in } \Gamma_\text{CYL} \cup \Gamma_\text{S} \cup \Gamma_\text{N} \times [0, T],\\
\end{cases}
\label{ns}
\end{equation}
where $\Omega$ is the 2-dimensional rectangular domain containing a circular hole (which represents the cylinder), depicted in Figure~\ref{fig:stokes_domain}. We set $T = 5000$ and $\Delta t = 0.2$ for time discretization in the full-order model.
The parameter $\mu$ controls the horizontal component of the velocity of the fluid entering the inlet (the vertical component is set to zero), and consequently also the Reynolds number. Since the parameter space is set to $\mathbb{P} = [0.0001, 0.002]$, we are here in laminar regime, making unnecessary the employment of turbulence models.
We impose on the physical walls (namely the cylinder and the two horizontal walls) the no-slip boundary condition.
For the training set, we considered 14 equispaced values of $\mu \in \mathbb{P}$ and $200$ time instants in the time window $[4000, 4800]$ with $\Delta t = 4$. Initial time steps show the propagation of boundary condition, obfuscating the real physical phenomena, and then they are discarded by the training set.

\begin{figure}[b!]
    \centering
    \includegraphics[width=.8\textwidth]{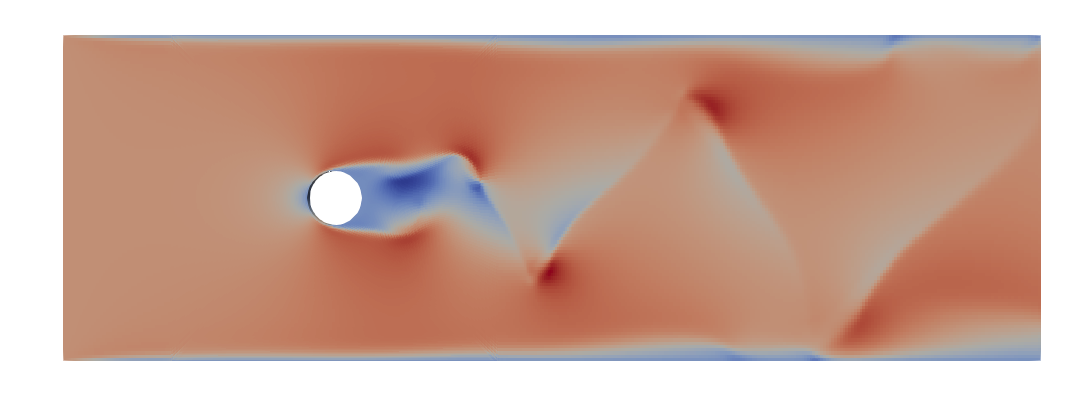}
    \label{cylinder:system}
    \caption{The system for $\mu=0.0015$ (one of the tested parameters) at $t=850$.}
\end{figure}

\noindent We employed the \emph{partitioned} procedure discussed in \cref{meth:partial} to test the results of our approach on the two untested parameters $\mu=0.0005$, $\mu=0.001775$, for which we predicted an approximated solution for the time instants between $4800$ and $5000$. For this experiment we only considered the horizontal component of the point-wise velocity of the fluid, in order to simplify the presentation of the results.

\begin{figure}[t!]
        \includegraphics[width=.47\textwidth]{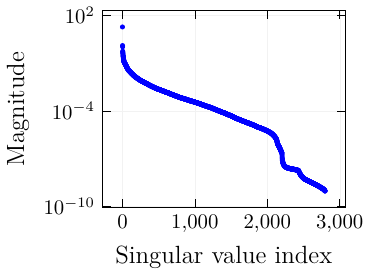}
        \includegraphics[width=.42\textwidth]{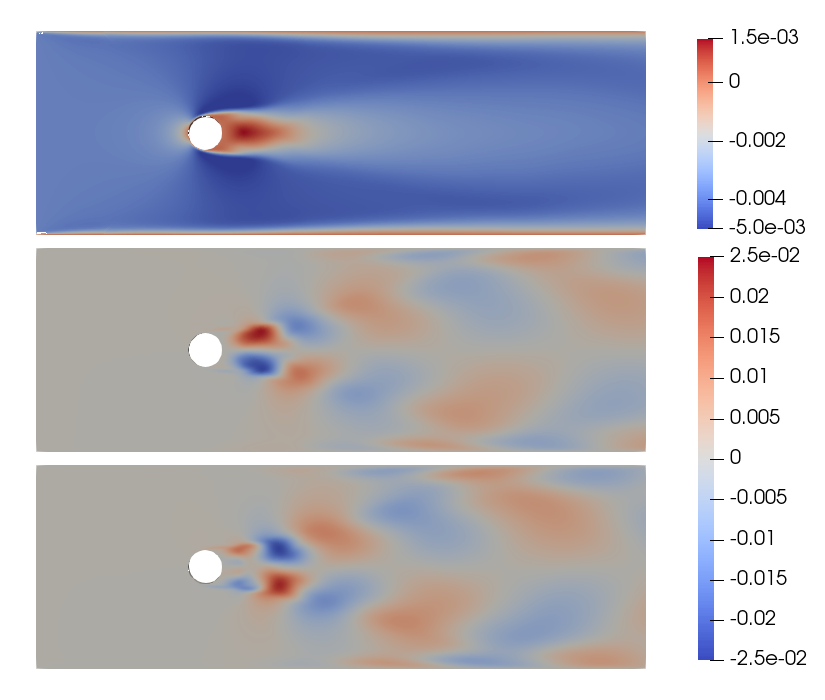}
    \caption{The singular values (left) and the first three POD modes (right) of the Navier-Stokes problem.}\label{fig:pod_ns}
\end{figure}

\begin{figure}[H]
\centering
\includegraphics[width=\textwidth]{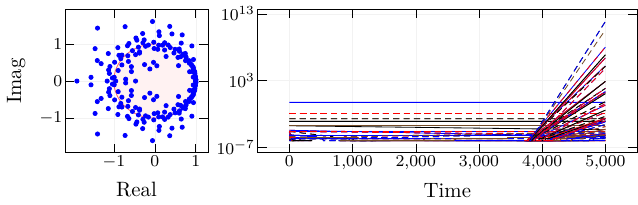}
\includegraphics[width=\textwidth]{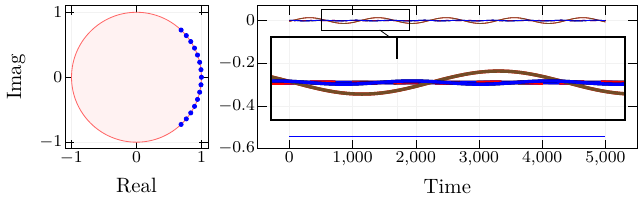}
\caption{DMD Stabilization in the Navier-Stokes problem. On top, DMD eigenvalues and dynamics before stabilization. On bottom, DMD eigenvalues and dynamics after stabilization.}\label{fig:stabil}
\end{figure}

We considered 40 POD modes, and did not perform an additional reduction during the application of DMD. We used HODMD, the variant of DMD described in \cref{sec:hodmd}, since the spatial complexity (which is 40 in this case, since we only retain 40 POD coefficients) was much lower than the number of time instants considered for the training (we remind to the reader that for the modular approach DMD is applied individually to each matrix $\matr{A}_{\mu_i} \in \mathbb{R}^{\npod \times \timedim}$ of POD coefficients).
We also used the stabilization approach presented in \cref{meth:stab} and retained only DMD modes corresponding to eigenvalues whose distance from the unit circle was at most $\epsilon = 10^{-3}$.
As we can see in \cref{fig:stabil}, thanks to such approach we are able to discard all the divergent modes, obtaining a good prediction also in test time window $[4800, 5000]$ (relative mean error equal to $0.04$). We remark that without stabilization the error for future time instants is several order of magnitude bigger, since the exponential growth of some modes.
\pgfplotsset{compat=1.5.1}

\begin{figure}[b!]
    \centering
    \includegraphics[width=\textwidth]{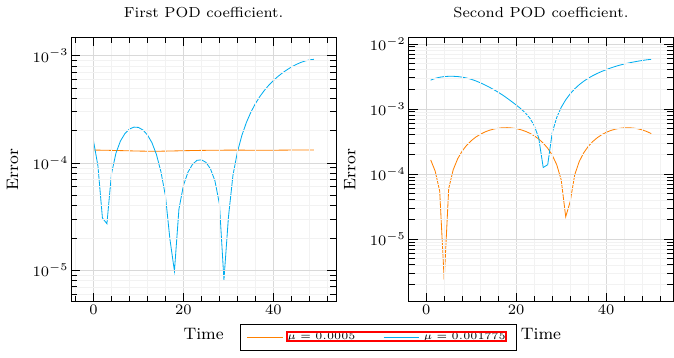}
    \caption{Absolute error for the first (on the left) and second (on the right) POD coefficients.}
    \label{cylinder:pod_error}
\end{figure}

\begin{figure}[t!]
    \includegraphics[width=\textwidth]{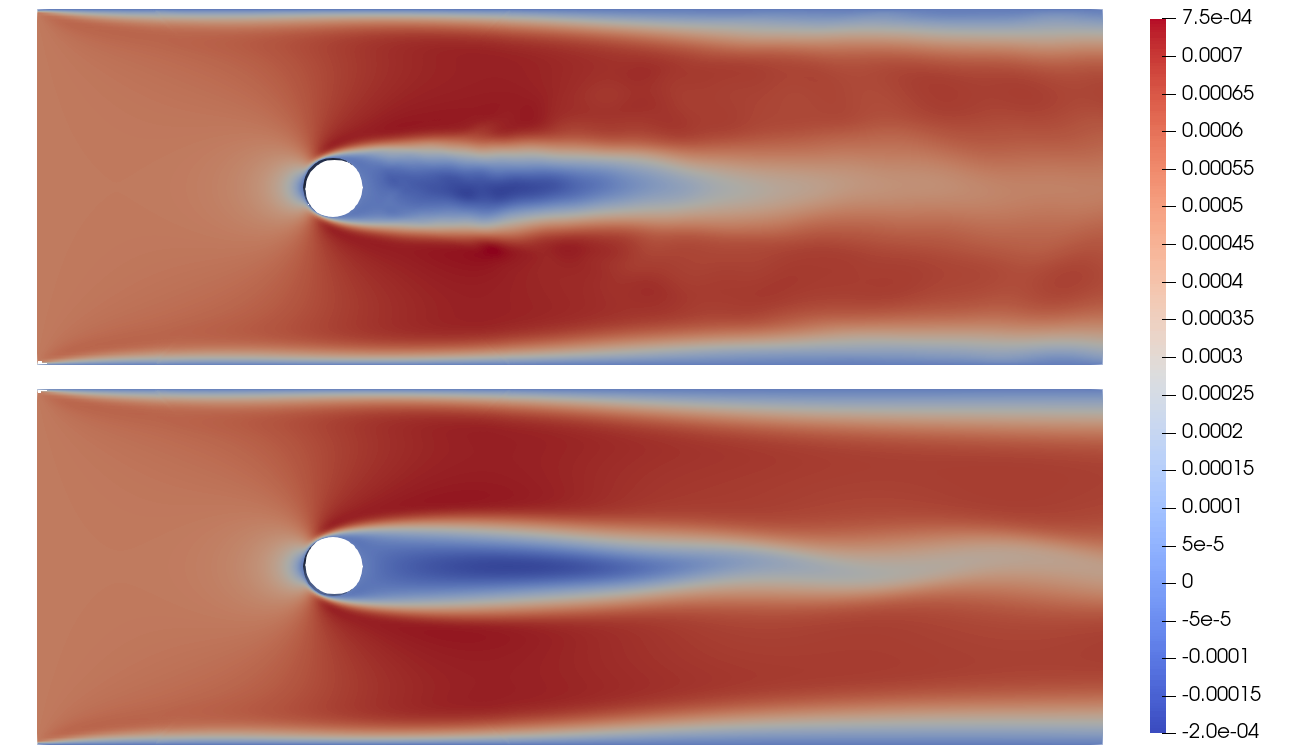}
    \caption{Results of the proposed algorithm (above) for the Cylinder problem for $\mu=0.0005$ at $t=1250$ (truth solution below).}
    \label{cylinder:untest_1_result}
\end{figure}

In \cref{cylinder:system} we show the physical appearance of the system for one of the tested parameters. After the examination of the decay of the singular values of the matrix $\matr{\mathcal{X}}_1$ (\cref{fig:pod_ns}) we observe that the system is poorly represented if we use a low number of POD modes. However, the proposed approach relies on the interpolation of POD modal coefficients, which becomes less and less precise (for same number of samples) at higher dimensions. Therefore we expect an increase in the error if the number of POD modes is increased. While we see some advantages in raising the number of POD modes  (especially for systems like this one) it is recommended to look for a compromise between the error given by regression and error given by POD space projection.

In \cref{fig:pod_ns} we also observe the first POD modes, which are the characteristic modes for low values of the parameter (first mode, in this case the system is stationary and has almost no waves) and high values of the parameter (second and third modes). This confirms that POD is able to extract dominant structures simultaneously from different instances of a parameterized system.

In \cref{cylinder:pod_error} we provide an analysis of the error on the first POD coefficients as a function of the time, in the predicted time frame ($4800 \leq t \leq 5000$). Like in the previous cases, the error is computed with a multiplication between the high-dimensional matrix which represents the system in the untested parameters and the matrix $\matr{U}_r$ of POD modes obtained during the \emph{offline phase}. Again, POD coefficients do not appear to be placed along a linear path as $\mu$ varies, however we are able to obtain an acceptable approximation for intermediate values of the parameter.

\begin{figure}[b!]
    \includegraphics[width=\textwidth]{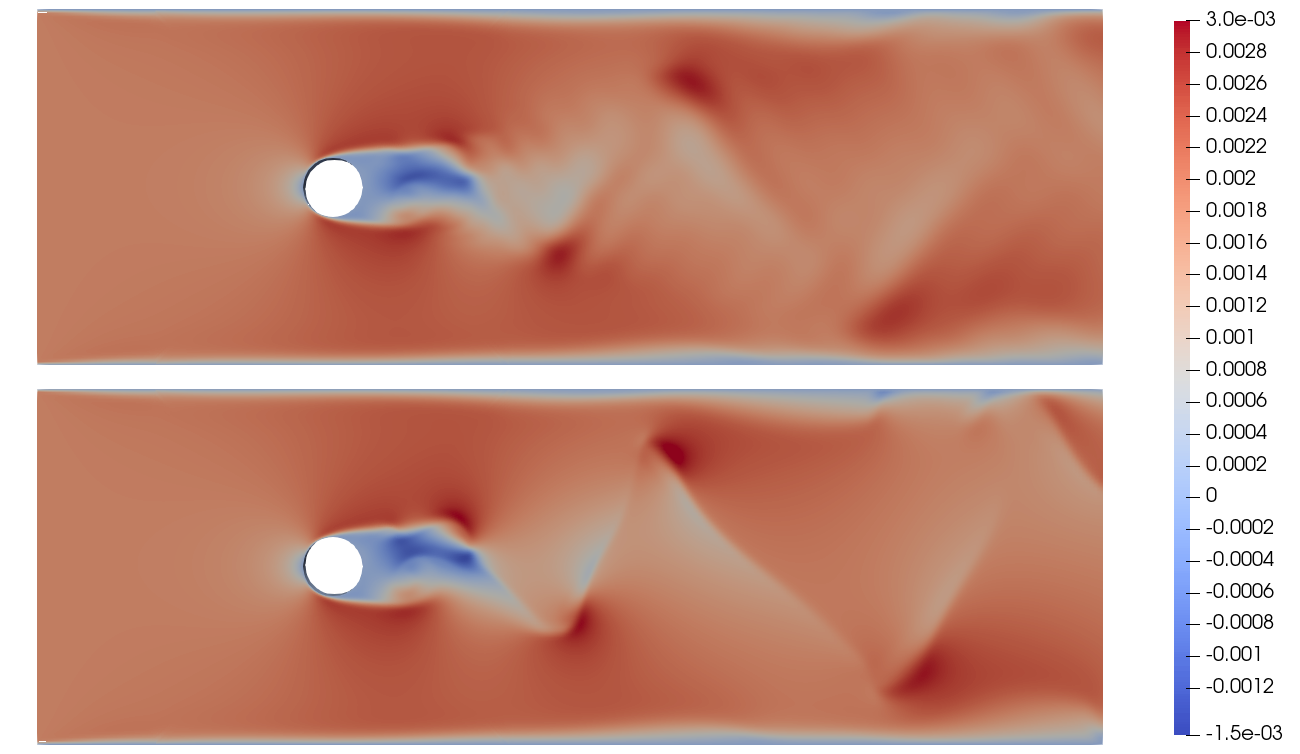}
    \caption{Results of the proposed algorithm (above) for the Cylinder problem for $\mu=0.001775$ at $t=5000$ (truth solution below).}
    \label{cylinder:untest_2_result}
\end{figure}%

Finally, in \cref{cylinder:untest_1_result,cylinder:untest_2_result} we show the results of the proposed algorithm (above) against the truth solution (below). In the regions outside the cylinder trail the pointwise error is higher, due to the low number of POD modes we considered. However we are able to obtain an acceptable approximation of the region of $\Omega$ in the trail: we obtain a low pointwise error near the cylinder; as we procede towards the right wall the error increases, but we are able to capture and represent the characteristic behaviors of the waves.

\section{Conclusions}
\label{sec:conclusions}
In this contribution, we introduced a novel approach to exploit dynamic mode decomposition even in the case of parametrized dynamical systems. Keeping the data-driven nature and the offline--online subdivision of the original algorithm, we exploit reduction and regression techniques in order to obtain a quick approximation of the output of interest for new untested parameters, in any temporal instant.
The such methodology demonstrated its effectiveness with three different numerical experiments, presenting some limitations where the original model is not well represented in a POD space, as for the Navier-Stokes problem. A possible solution for these problems has been proposed in~\cite{martin2021}, in parallel to the preparation of this manuscript, by interpolating the different DMD operators computed for different parameter values.

\RB{
\label{conc:future}
In this work, we presented the results obtained by taking parameters individually (\emph{partitioned} approach) and all together (\emph{monolithic} approach). We theorize that the accuracy may be improved by introducing some clustering techniques to group parameters (e.g. spatial position, resulting Reynolds number, \dots). The goal is to put into the same cluster parameters which lead to similar realizations of the dynamical system. This modification -- which could be tuned according to some hyperparameter that determines the number of clusters to be found in the set -- might enable finding the optimal approach by moving on a segment whose extremes are the \emph{monolithic} and the \emph{partitioned} approaches. It may even be interesting to consider clusters with non-empty intersections, in order to improve the linking between different groups of parameters. However, we expect that finding an optimal configuration, in this case, would be much more difficult. Another further improvement of the proposed framework could be the usage of other dimensionality reduction techniques in order to replace the POD method. A possible candidate is a class of \emph{autoencoder} methods, which have already been used to analyze and explain the temporal evolution of fluids \cite{agostini2020exploration, omata2019novel, fukami2020convolutional}. Finally, we postponed the application of noise-tolerant DMD variants in order to improve the robustness and accuracy of our method.
}

\section*{Acknowledgements}
This work was partially supported by European Union Funding for Research and Innovation ---
Horizon 2020 Program --- in the framework of European Research Council
Executive Agency: H2020 ERC CoG 2015 AROMA-CFD project 681447 ``Advanced
Reduced Order Methods with Applications in Computational Fluid Dynamics'' P.I.
Gianluigi Rozza.

\bibliographystyle{abbrv}
\bibliography{biblio}

\begin{thebibliography}{10}

\bibitem{agostini2020exploration}
L.~Agostini.
\newblock Exploration and prediction of fluid dynamical systems using
  auto-encoder technology.
\newblock {\em Physics of Fluids}, 32(6):067103, 2020.

\bibitem{farhat}
D.~Amsallem and C.~Farhat.
\newblock Interpolation method for adapting reduced-order models and
  application to aeroelasticity.
\newblock {\em AIAA journal}, 46(7):1803--1813, 2008.

\bibitem{amsallem2011online}
D.~Amsallem and C.~Farhat.
\newblock An online method for interpolating linear parametric reduced-order
  models.
\newblock {\em SIAM Journal on Scientific Computing}, 33(5):2169--2198, 2011.

\bibitem{askham2018variable}
T.~Askham and J.~N. Kutz.
\newblock Variable projection methods for an optimized dynamic mode
  decomposition.
\newblock {\em SIAM Journal on Applied Dynamical Systems}, 17(1):380--416,
  2018.

\bibitem{bagheri2013effects}
S.~Bagheri.
\newblock Effects of small noise on the dmd/koopman spectrum.
\newblock {\em Bulletin Am. Phys. Soc}, 58(18):H35, 2013.

\bibitem{nearest}
N.~Bhatia et~al.
\newblock Survey of nearest neighbor techniques.
\newblock {\em arXiv preprint arXiv:1007.0085}, 2010.

\bibitem{brunton2014compressive}
S.~L. Brunton, J.~H. Tu, I.~Bright, and J.~N. Kutz.
\newblock Compressive sensing and low-rank libraries for classification of
  bifurcation regimes in nonlinear dynamical systems.
\newblock {\em SIAM Journal on Applied Dynamical Systems}, 13(4):1716--1732,
  2014.

\bibitem{dawson2016characterizing}
S.~Dawson, M.~S. Hemati, M.~O. Williams, and C.~W. Rowley.
\newblock Characterizing and correcting for the effect of sensor noise in the
  dynamic mode decomposition.
\newblock {\em Experiments in Fluids}, 57(3):1--19, 2016.

\bibitem{DemoOrtaliGustinRozzaLavini2020BUMI}
N.~Demo, G.~Ortali, G.~Gustin, G.~Rozza, and G.~Lavini.
\newblock An efficient computational framework for naval shape design and
  optimization problems by means of data-driven reduced order modeling
  techniques.
\newblock {\em Bollettino dell'Unione Matematica Italiana}, Nov 2020.

\bibitem{demo2018shape}
N.~Demo, M.~Tezzele, G.~Gustin, G.~Lavini, and G.~Rozza.
\newblock Shape optimization by means of proper orthogonal decomposition and
  dynamic mode decomposition. submitted.
\newblock {\em arXiv preprint arXiv:1803.07368}, 2018.

\bibitem{pydmd}
N.~Demo, M.~Tezzele, and G.~Rozza.
\newblock Pydmd: Python dynamic mode decomposition.
\newblock {\em Journal of Open Source Software}, 3(22):530, 2018.

\bibitem{fukami2020convolutional}
K.~Fukami, T.~Nakamura, and K.~Fukagata.
\newblock Convolutional neural network based hierarchical autoencoder for
  nonlinear mode decomposition of fluid field data.
\newblock {\em Physics of Fluids}, 32(9):095110, 2020.

\bibitem{gadalla2021les}
M.~Gadalla, M.~Cianferra, M.~Tezzele, G.~Stabile, A.~Mola, and G.~Rozza.
\newblock {On the comparison of LES data-driven reduced order approaches for
  hydroacoustic analysis}.
\newblock {\em {Computers \& Fluids}}, 216:104819, 2021.

\bibitem{georgaka2018parametric}
S.~Georgaka, G.~Stabile, G.~Rozza, and M.~J. Bluck.
\newblock Parametric pod-galerkin model order reduction for unsteady-state heat
  transfer problems.
\newblock {\em arXiv preprint arXiv:1808.05175}, 2018.

\bibitem{grepl2007efficient}
M.~A. Grepl, Y.~Maday, N.~C. Nguyen, and A.~T. Patera.
\newblock Efficient reduced-basis treatment of nonaffine and nonlinear partial
  differential equations.
\newblock {\em ESAIM: Mathematical Modelling and Numerical Analysis},
  41(3):575--605, 2007.

\bibitem{hemati2017biasing}
M.~S. Hemati, C.~W. Rowley, E.~A. Deem, and L.~N. Cattafesta.
\newblock De-biasing the dynamic mode decomposition for applied koopman
  spectral analysis of noisy datasets.
\newblock {\em Theoretical and Computational Fluid Dynamics}, 31(4):349--368,
  2017.

\bibitem{martin2021}
M.~W. Hess, A.~Quaini, and G.~Rozza.
\newblock A data-driven surrogate modeling approach for time-dependent
  incompressible navier-stokes equations with dynamic mode decomposition and
  manifold interpolation.
\newblock {\em Submitted}, 2021.

\bibitem{rozza_pod}
J.~S. Hesthaven, G.~Rozza, B.~Stamm, et~al.
\newblock {\em Certified reduced basis methods for parametrized partial
  differential equations}, volume 590.
\newblock Springer, 2016.

\bibitem{hijazi2020data}
S.~Hijazi, G.~Stabile, A.~Mola, and G.~Rozza.
\newblock Data-driven pod-galerkin reduced order model for turbulent flows.
\newblock {\em Journal of Computational Physics}, 416:109513, 2020.

\bibitem{hoang2021projection}
C.~Hoang, K.~Chowdhary, K.~Lee, and J.~Ray.
\newblock Projection-based model reduction of dynamical systems using
  space-time subspace and machine learning.
\newblock {\em arXiv preprint arXiv:2102.03505}, 2021.

\bibitem{koopman}
B.~O. Koopman.
\newblock Hamiltonian systems and transformation in hilbert space.
\newblock {\em Proceedings of the national academy of sciences of the united
  states of america}, 17(5):315, 1931.

\bibitem{kramer2017sparse}
B.~Kramer, P.~Grover, P.~Boufounos, S.~Nabi, and M.~Benosman.
\newblock Sparse sensing and dmd-based identification of flow regimes and
  bifurcations in complex flows.
\newblock {\em SIAM Journal on Applied Dynamical Systems}, 16(2):1164--1196,
  2017.

\bibitem{dmd}
J.~N. Kutz, S.~L. Brunton, B.~W. Brunton, and J.~L. Proctor.
\newblock {\em Dynamic mode decomposition: data-driven modeling of complex
  systems}.
\newblock SIAM, 2016.

\bibitem{hodmd}
S.~Le~Clainche and J.~M. Vega.
\newblock Higher order dynamic mode decomposition.
\newblock {\em SIAM Journal on Applied Dynamical Systems}, 16(2):882--925,
  2017.

\bibitem{hodmdflow}
S.~Le~Clainche and J.~M. Vega.
\newblock Higher order dynamic mode decomposition to identify and extrapolate
  flow patterns.
\newblock {\em Physics of Fluids}, 29(8):084102, 2017.

\bibitem{hodmdnoisy}
S.~Le~Clainche, J.~M. Vega, and J.~Soria.
\newblock Higher order dynamic mode decomposition of noisy experimental data:
  The flow structure of a zero-net-mass-flux jet.
\newblock {\em Experimental Thermal and Fluid Science}, 88:336--353, 2017.

\bibitem{lu2020prediction}
H.~Lu and D.~M. Tartakovsky.
\newblock Prediction accuracy of dynamic mode decomposition.
\newblock {\em SIAM Journal on Scientific Computing}, 42(3):A1639--A1662, 2020.

\bibitem{omata2019novel}
N.~Omata and S.~Shirayama.
\newblock A novel method of low-dimensional representation for temporal
  behavior of flow fields using deep autoencoder.
\newblock {\em Aip Advances}, 9(1):015006, 2019.

\bibitem{penrose1956best}
R.~Penrose.
\newblock On best approximate solutions of linear matrix equations.
\newblock In {\em Mathematical Proceedings of the Cambridge Philosophical
  Society}, volume~52, pages 17--19. Cambridge University Press, 1956.

\bibitem{quarteroni}
A.~Quarteroni, R.~Sacco, and F.~Saleri.
\newblock {\em Numerical mathematics}, volume~37.
\newblock Springer Science \& Business Media, 2010.

\bibitem{RozzaMalikDemoTezzeleGirfoglioStabileMola2018ECCOMAS}
G.~Rozza, M.~H. Malik, N.~Demo, M.~Tezzele, M.~Girfoglio, G.~Stabile, and
  A.~Mola.
\newblock {Advances in Reduced Order Methods for Parametric Industrial Problems
  in Computational Fluid Dynamics}.
\newblock In {\em Proceedings of the ECCOMAS Congress 2018}. ECCOMAS, ECCOMAS,
  2018.

\bibitem{sashidhar2022bagging}
D.~Sashidhar and J.~N. Kutz.
\newblock Bagging, optimized dynamic mode decomposition for robust, stable
  forecasting with spatial and temporal uncertainty quantification.
\newblock {\em Philosophical Transactions of the Royal Society A},
  380(2229):20210199, 2022.

\bibitem{pdmd}
T.~Sayadi, P.~J. Schmid, F.~Richecoeur, and D.~Durox.
\newblock Parametrized data-driven decomposition for bifurcation analysis, with
  application to thermo-acoustically unstable systems.
\newblock {\em Physics of Fluids}, 27(3):037102, 2015.

\bibitem{TezzeleDemoMolaRozza2018}
M.~Tezzele, N.~Demo, A.~Mola, and G.~Rozza.
\newblock {\em An integrated data-driven computational pipeline with model
  order reduction for industrial and applied mathematics}.
\newblock 2018.

\bibitem{TezzeleDemoStabileMolaRozza2020MOR}
M.~Tezzele, N.~Demo, G.~Stabile, A.~Mola, and G.~Rozza.
\newblock Enhancing cfd predictions in shape design problems by model and
  parameter space reduction.
\newblock {\em Advanced Modeling and Simulation in Engineering Sciences},
  7(1):40, Oct 2020.

\bibitem{vasconcelos2019dynamic}
E.~Vasconcelos~Filho and P.~L. dos Santos.
\newblock A dynamic mode decomposition approach with hankel blocks to forecast
  multi-channel temporal series.
\newblock {\em IEEE Control Systems Letters}, 3(3):739--744, 2019.

\bibitem{gaussian_process}
C.~K. Williams and C.~E. Rasmussen.
\newblock {\em Gaussian processes for machine learning}, volume~2.
\newblock MIT press Cambridge, MA, 2006.

\end{thebibliography}

\end{document}